\documentclass[12pt,notitlepage,twoside]{article}
\usepackage{latexsym}
\usepackage{amsmath}
\usepackage{amsfonts}
\usepackage{amssymb}

\renewcommand{\a }{\alpha }
\renewcommand{\b }{\beta }

\newcommand{\f}{\varphi }

\newcommand{\ii}{{\rm i}}

\newcommand{\e }{\varepsilon }

\newcommand{\x }{\xi}

\renewcommand{\l }{\lambda }

\renewcommand{\r }{\rho }

\renewcommand{\O }{\Omega }

\newcommand{\ov}{\overline}
\newcommand{\wtilde }{\widetilde}

\newcommand{\be}{\begin{equation}}
\newcommand{\ee}{\end{equation}}

\newcommand{\torus}{\mathbb{T}}
 
\newcommand{\hv}{{\hat v}}
\newcommand{\hf}{{\hat \f}}

\newcommand{\intt}{\int_0^{2\pi}}
\newcommand{\intOa}{\int_{\O_\a}}

\newcommand{\tw}{\widetilde{w}}

\newcommand{\tf}{\widetilde{f}}
\newcommand{\bv}{\bar{v}}
\newcommand{\hbv}{\hat{\bar{v}}}
\newcommand{\bw}{\bar{w}}
\newcommand{\bu}{\bar{u}}

\newcommand{\uc}{\underline{c}}

\newcommand{\oO}{\overline{\O}}
\newcommand{\qed}{\hfill$\Box$}

\newtheorem{thm}{Theorem}[section]
\newtheorem{pro}[thm]{Proposition}
\newtheorem{lem}[thm]{Lemma}
\newtheorem{rem}[thm]{Remark}

\newtheorem{df}[thm]{Definition}
\numberwithin{equation}{section}

\newcommand\equ[1]{{\rm (\ref{#1})}}

\author{{\sc  Massimiliano Berti, Luca Biasco}
\thanks{ Massimiliano Berti, SISSA, via Beirut 2-4,
Trieste, Italy, \texttt{berti@sissa.it} \ . Luca Biasco,
Universit\`a di Roma $3$, Largo S. Leonardo Murialdo, Roma, Italy,
\texttt{biasco@mat.uniroma3.it} \ .}}

\title { \Large \textbf{Forced vibrations of
wave equations with non-monotone nonlinearities}}

\begin{document}

\date{}

\maketitle

{\footnotesize
\begin{abstract}

\noindent
We prove existence and regularity of periodic in time solutions of
completely resonant nonlinear forced wave equations with Dirichlet
boundary conditions for a large class of non-monotone
forcing terms.
Our approach is based on
a variational Lyapunov-Schmidt reduction.
It turns out that the infinite dimensional
bifurcation equation exhibits an intrinsic
lack of compactness.
We solve it via a minimization argument
and a-priori estimate methods inspired to regularity theory of
[R67].


\end{abstract}
}

\tableofcontents

\section{Introduction}\label{sec:intro}

In this paper we consider the problem of finding
nontrivial time-periodic solutions of
the completely resonant nonlinear {\sl forced} wave equation
\begin{equation}\label{equ}
\square u = \e f ( t, x, u; \e )
\end{equation}
with Dirichlet boundary conditions
\begin{equation}\label{bc}
u ( t, 0 ) = u ( t, \pi ) = 0
\end{equation}
where $ \square := \partial_{tt} - \partial_{xx} $
is the D'Alembertian operator,
$ \e $ is a small parameter and
the nonlinear forcing  term $ f (t,x,u; \e ) $ is $T$-periodic in time.
We consider the case when
$ T $
is a rational multiple of $ 2 \pi $
and, for simplicity of exposition,
we shall assume
$$
T = 2 \pi \ .
$$
We look for nontrivial $ 2\pi $-periodic in time solutions $u(t,x) $
of \equ{equ}-\equ{bc}, i.e. satisfying
\begin{equation}\label{pc}
u( t + 2 \pi, x ) = u(t,x) \ .
\end{equation}

For $\e = 0 $, \equ{equ}-\equ{bc} reduces to the linear
homogeneous wave
equation
\begin{equation}\label{equlin}
\begin{cases}
\square u = 0 \cr
u ( t, 0 ) = u ( t, \pi ) = 0
\end{cases}
\end{equation}
which possesses an {\sl infinite} dimensional space of solutions
which are $ 2 \pi $-periodic in time and of the form
$ v(t,x)= $ $ \hv(t+x)- $ $ \hv(t-x) $ for any $2 \pi$-periodic function
$ \hv ( \cdot) $.
For this reason equation \equ{equ}-\equ{bc} is called
completely resonant.
\\[1mm]
\indent
The main difficulty for proving existence of solutions of
\equ{equ}-\equ{bc}-\equ{pc}
for $ \e \neq 0 $
is to find from which periodic orbits
of the linear equation
\equ{equlin} the
 solutions of the nonlinear equation
\equ{equ}
branch off. This requires to solve an infinite dimensional
bifurcation equation (also called  kernel equation)
with an intrinsic {\sl lack of
compactness}.
\\[1mm]
\indent
The first breakthrough regarding problem \equ{equ}-\equ{bc}-\equ{pc}
was achieved
by Rabinowitz in [R67] where existence and regularity of solutions
was proved for nonlinearities satisfying the strongly monotone
assumption  $ (\partial_u f) (t,x,u) \geq \b > 0 $. Using methods
inspired by the theory of elliptic regularity,
[R67]  proved
the existence of a unique curve of smooth solutions for $ \e $ small.
Other existence
results of weak and classical solutions
have been obtained,
still in the strongly monotone case,
in [DST68]-[L69]-[BN78].

Subsequently, Rabinowitz [R71] was able to prove
existence of weak solutions of \equ{equ}-\equ{bc}-\equ{pc}
for a class of weakly monotone nonlinearities
like $ f(t,x,u) = u^{2k+1} + G(t,x,u) $ where
$ G (t,x,u_2) \geq G(t,x,u_1) $ if $ u_2 \geq u_1 $.
Actually, in [R71] bifurcation of a global continuum branch of
weak solutions is proved.
For other local existence results in the weakly monotone case
we mention [T69]-[H70].
\\[1mm]
\indent
In all the quoted papers
the monotonicity assumption (strong or weak) is
the key property
for overcoming the
lack of compactness in
the infinite dimensional kernel equation.
\\[1mm]
\indent
We underline that, in general, the weak solutions obtained
in [R71] are only continuous functions.
Concerning
regularity,
Brezis and Nirenberg [BN78] proved -but only for strongly
monotone nonlinearities-
that any $ L^\infty $-solution of \equ{equ}-\equ{bc}-\equ{pc}
is smooth, even in the nonperturbative case
$ \e = 1 $,  whenever
the  nonlinearity
 $ f $ is smooth.
\\[1mm]
\indent
On the other hand,
very little is known about existence and regularity of solutions
if we drop the monotonicity assumption on the forcing term $ f $.
Willem [W81], Hofer [H82] and Coron [C83] have considered
the class of equations \equ{equ}-\equ{bc}
where $ f (t,x,u) = g(u) + h (t,x) $, $ \e = 1 $, and
$ g(u) $ satisfies suitable linear growth conditions.
Existence of weak solutions is proved, in [W81]-[H82],
for a set of $ h $ dense in $ L^2 $,
although explicit criteria that characterize such $ h $
are not provided.
The infinite dimensional bifurcation
problem is overcome by assuming non-resonance hypothesys
between the asymptotic
behaviour of $ g(u) $ and the spectrum of $ \square $.
On the other side,
Coron [C83] finds weak solutions
assuming the  additional symmetry $ h(t,x)= h ( t + \pi , \pi - x )$
and restricting to the space of functions satisfying
$ u(t,x)= u ( t + \pi , \pi - x ) $, where
the Kernel of the d'Alembertian operator
$ \square $ reduces to $ 0 $.
For some more recent results see for example
[BDL99].
\\[1mm]
\indent
In the present paper we prove existence
and regularity of solutions of \equ{equ}-\equ{bc}-\equ{pc}
for a large class of
{\sc nonmonotone} forcing terms $ f(t,x,u) $, including, for example,
\begin{description}
\item $f(t,x,u) = \pm u^{2k} +  h(t,x) $,
see {\it Theorem 1} \ ;
\item $  f(t,x,u) =$ $ \pm u^{2k} + $ $u^{2k+1} + $ $ h(t,x) $,
see {\it Theorem 2} \ ;
\item $ f (t,x,u) = \pm u^{2k} + {\wtilde f}(t,x,u) $  \ with \
$ (\partial_u {\wtilde f})(t,x,u) \geq \b > 0 $,
see {\it  Theorem 3} \ .
\end{description}

The precise results will be stated in the next
subsection \ref{subsec:main}, see {\it
Theorems 1, 2} and {\it  3}.
Their proof
is based on
a variational Lyapunov-Schmidt reduction, minimization
arguments and a-priori estimate methods inspired
to regularity theory of [R67].
We anticipate that
our approach -explained in subsection \ref{vlsr}-
is not
merely a sharpening of the ideas of [R67]-[R71],
which, to deal with non monotone nonlinearities,
require a significant
change of prospective.
\\[1mm]
\indent
We mention that in the last years several results on
bifurcation of free vibrations for completely resonant
autonomous wave equations
have been proved in [B99]-[BP01]-[BB03]-[BB04a]-[BB04b]-[GMP04].
The main differences with respect to the present case
are that: a ``small divisor'' problem in solving
the ``range equation'' appears (here no
small divisor problem is present
due to the assumption $ T = 2 \pi$, see remark \ref{orto}), but
the infinite dimensional
``bifurcation equation'' -whose solutions
is the main problem of the present
paper- gains crucial compactness
properties, see remark \ref{rem:aut}.

\subsection{Main Results}\label{subsec:main}

We look for solutions $ u : \O \to \mathbb{R} $ of \equ{equ}
in the Banach space
$$
E := H^1(\O ) \cap C^{1\slash 2}_0 (\oO ), \qquad \O :=
\mathbb{T} \times (0,\pi)
$$
where $H^1(\O ) $ is the usual Sobolev space
and $ C^{1/2}_0 (\oO) $ is
the space of all the $1/2$-H\"{o}lder continuous
functions $u:\oO\longrightarrow \mathbb{R}$
satisfying \equ{bc},
endowed with norm\footnote{Here $ \|u\|_{H^1( \O ) }^2 :=$
$\| u \|_{L^2(\O)}^2 + $ $\|u_x\|_{L^2(\O)}^2 + $ $ \|u_t\|_{L^2(\O)}^2 $
and
$$
\|u\|_{C^{1 \slash 2 }( \oO)} := \| u \|_{ C^0 ( \O )} +
\sup_{(t,x) \neq (t_1,x_1)} \frac{| u(t,x) - u(t_1,x_1) |}{(|
t - t_1| + | x - x_1|)^{1\slash 2} \label{hn}   } \ .
$$}
$$
\| u \|_E := \|u\|_{H^1(\O)} + \|u\|_{C^{1 \slash 2} (\oO)} \ .
$$
\noindent
Critical points of the Lagrangian action functional
$ \Psi \in C^1(E,\mathbb{R})$
\begin{equation}\label{Psi}
\Psi(u):=\Psi(u,\e):=\int_\O
\Big[ \frac{u_t^2}{2}-\frac{u_x^2}{2}+\e F(t,x,u; \e ) \Big] \ dt dx \ ,
\end{equation}
where $F( t,x,u; \e ):=\int_0^u f(t,x,\xi; \e)d\xi $,
are weak solutions of \equ{equ}-\equ{bc}-\equ{pc}.
\\[1mm]
\indent
For $ \e = 0 $,  the critical points of
$\Psi $ in $E$ reduce to the solutions of the linear equation
\equ{equlin} and form the
subspace $V := N \cap H^1 (\O) $
where
\begin{equation}\label{KerN}
N :=\left\{\
v(t,x)= \hv(t+x)-\hv(t-x) 
\  \Big|\ \ \hv\in L^2(\mathbb{T}) \ {\rm and} \
\int_0^{2\pi} \hv (s ) \ ds =0\
  \right\}.
\end{equation}
Note that $ V := $ $ N \cap H^1 (\O) =$
$ \{ v(t,x) = $ $ \hv(t+x)-\hv(t-x) \in N $ $ | $ $
\hv \in H^1(\mathbb{T})\} \subset E $, since
any function $ \hv \in H^1(\mathbb{T})$ is
$ 1 \slash 2$-H\"older continuous.
\\[1mm]
Let $ N^\bot := \{ h \in L^2 (\O) \ | \ \int_\O h v = 0, \ \forall v \in N
\} $ denote the $L^2 (\O) $-orthogonal of $ N $.
\\[1mm]
We prove the following Theorem:

\medskip
\noindent
{\bf Theorem 1}
{\it Let
$f(t,x,u)= \b u^{2k} + h(t,x) $ 
and
$ h \in N^\bot 
$ satisfies
$ h(t,x)> 0 $ $($or $ h(t,x) <0 $$)$ a.e. in $ \O $.
Then, for $ \e $ small enough, there exists at least one weak solution
$ u \in E $
of \equ{equ}-\equ{bc}-\equ{pc} with $\|u\|_E  
\leq C | \e |$. If, moreover,
$ h \in H^j ( \O) \cap C^{j-1} ( \oO ) $,
$ j \geq 1 $, then $ u \in H^{j+1} (\O) \cap C^j_0 (\oO)$
with
$\| u\|_{H^{j+1}(\O)}+\| u\|_{C^j(\oO)}\leq C |\e|$
and therefore, for $j\geq 2,$ $ u $ is a classical solution.}

\medskip

Theorem 1 is a Corollary of the
following more general result which enables to deal
with non-monotone nonlinearities like, for example,
$ f(t,x,u)= \pm (\sin x ) \ u^{2k} + h (t,x)$,
$ f(t,x,u)= \pm  \ u^{2k} + u^{2k+1} + h (t,x) $.

\medskip

\noindent
{\bf Theorem 2}
{\it Let $f(t,x,u)= g(t,x,u) + h(t,x)$,
 $ h (t,x) \in N^\bot $ and 
$$
g(t,x,u):= \b(x) u^{2k} +\mathcal{R}(t,x,u)
$$
where
$\mathcal{R},$ $\partial_t\mathcal{R}$, $\partial_u\mathcal{R} \in
C (\oO \times \mathbb{R}, \mathbb{R} )$
satisfy\footnote{The notation
$ f(z) = o(z^p )$, $ p \in \mathbb{N} $,
means that $ f(z) \slash |z|^p \to 0 $ as $z \to 0 $.
$f(z) = O(z^p )$ means that there exists a constant $C>0$
such that $|f(z)| \leq C|z|^p$ for all $z$ in a neighboorhood of $0$.}
\begin{equation}\label{h}
     \|\mathcal{R}(\cdot,u)\|_{C(\oO)}
     =o(u^{2k})\,,\ \
     \|\partial_t\mathcal{R}(\cdot,u)\|_{C(\oO)}
=O(u^{2k})\,,\ \
    \|\partial_u\mathcal{R}(\cdot,u)\|_{C(\oO)}
    =o(u^{2k-1})\, ,
\end{equation}
and $\b \in C([0,\pi], \mathbb{R})$ verifies, for $x \in (0,\pi )$,
$\b(x) > 0$ $($or $ \b(x) < 0 $$)$ and $ \b(\pi- x) = \b (x)$.
\\[1.5mm]
${{\bf (i)}}$ $(${\bf Existence}$)$ Assume there exists a
weak solution $ H \in E $ of $ \square  H = h $ such that
\begin{equation}\label{G:positiva}
    H(t,x) > 0
    \qquad  \ \mbox{$($or $ H(t,x) < 0 $$)$}  \  \qquad \forall
(t,x) \in  \O\ .
\end{equation}
Then, for $\e $ small enough,
there exists at least one weak solution $ u\in E $ of
\equ{equ}-\equ{bc}-\equ{pc} satisfying
$\| u \|_E \leq C |\e | $.
\\[1.5mm]
${{\bf (ii)}}$ $(${\bf Regularity}$)$  If, moreover,
$ h \in H^j ( \O )\cap C^{j-1}(\oO) $,
$\b\in H^j\big((0,\pi)\big),$
$\mathcal{R},$ $\partial_t\mathcal{R}$, $\partial_u\mathcal{R}$
$\in C^j (\oO\times \mathbb{R})$, $ j \geq 1 $,
then $ u \in H^{j+1} (\O) \cap C^j_0 (\oO)$
and, for  $ j \geq 2 $, $u$ is a classical solution.}

\medskip

Note that Theorem 2
does not require any  growth condition on $g$  at infinity. 
In particular it applies for any analytic
function $ g(u) $
satisfying $ g(0)= g'(0 ) = \ldots = g^{2k-1} (0) = 0 $ and
$g^{2k}(0) \neq 0 $.
\\[1mm]
\indent
We now collect some comments on the previous results.

\begin{rem}\label{rem:Hp1}
The assumption $h \in N^\bot$ is not of technical nature
both in Theorem 1 and in Theorem 2
$($at least if $g=g(x,u) = g(x,-u)= g(\pi-x,u) )$. Indeed, if
$ h \notin N^\bot $, periodic solutions of problem
\equ{equ}-\equ{bc}-\equ{pc} do not exist
in any fixed ball $\{ \| u \|_{L^\infty} \leq R \},$ $R>0$,
 for $\e $ small; see remark
\ref{rem:non}.
\end{rem}

\begin{rem}\label{rem:Hp2}
In Theorem 2
hypothesys \equ{G:positiva} and $ \beta > 0 $  $($or $\beta < 0 )$
are assumed to prove the existence
of a minimum of the ``reduced action functional'' $\Phi$, see
\equ{Phi}.
A sufficient condition implying
\equ{G:positiva} is
$ h > 0 $ a.e. in $ \O $,
see the ``maximum principle'' Proposition \ref{lem:6}.
This is also
the key step to derive Theorem 1 from Theorem 2.

\noindent
We also note 
that hypothesys \equ{G:positiva} can be weakened, see remark \ref{ipdeb}.
\end{rem}

\begin{rem}\label{rem:Hp3}
{\bf (Regularity)}
It is quite surprising that the weak solution $u$ of Theorems 1,
2 is actually smooth.
Indeed, while regularity always holds true
for strictly monotone nonlinearities (see [R67]-[BN78]),
yet for weakly monotone $ f $
it is not proved in general, unless the weak solution $u $ verifies
$\|\Pi_N u \|_{L^2} \geq C  > 0  $ (see [R71]). Note,
on the contrary,
that the weak solution $u$ of Theorem 2 satisfies
$\| \Pi_N u \|_{L^2} = O( \e ) $.

Moreover, assuming
\begin{equation}\label{hbis}
    \| \partial_t^l \partial_x^m \partial_u^n \mathcal{R}\|_{C(\oO)}
    = O(u^{2k-n})\,, \quad
    \forall\,
    0\leq l,n\leq j+1\,, 0\leq m\leq j\,,
    l+m+n\leq j+1
\end{equation}
we can also prove the estimate (see remark \ref{rem:inde})
\begin{equation}\label{sti:regol}
\| u\|_{H^{j+1}(\O)}+\| u\|_{C^j(\oO)}\leq C |\e|\,.
\end{equation}
\end{rem}

\begin{rem}\label{rem:Hp4}
{\bf (Multiplicity)}
For nonmonotone nonlinearities $f$ one can {\sc not}
in general expect unicity of the solutions.
Actually, for $ f(t,x,u) = g(x,u) + h(t,x) $ with
$ g(x,u) = g(x, -u) $, $ g (\pi - x, u) = g (x, u) $,
there exist infinitely many $ h \in N^\bot $
for which  problem \equ{equ}-\equ{bc}-\equ{pc} has (at least)
$3$ solutions, see remark \ref{multi}.
\end{rem}

Finally, we  extend
the result of [R67]
proving existence
of periodic solutions for nonmonotone nonlinearities $ f ( t, x, u ) $
obtained adding
to a
nonlinearity $\wtilde f (t,x,u)$ as in
[R67] (i.e. $ \partial_u \tf \geq \b > 0 $)
any nonmonotone term $ a(x,u) $ satisfying
\begin{equation}\label{a1}
a(x,-u)=a(x,u)\,,\qquad
a(\pi-x,u)=a(x,u)
\end{equation}
or
\begin{equation}\label{a2}
a(x,-u)=-a(x,u)\,,\qquad
a(\pi-x,u)=-a(x,u)\, .
\end{equation}

A prototype nonlinearity
is $ f(t,x,u)= \pm u^{2k} + \tf (t,x,u) $ with
$ \partial_u \tf \geq \b > 0 $.

\medskip

\noindent
{\bf Theorem 3}
{\it Let $ f ( t, x, u ) = $ $ \tf (t,x,u) + $ $ a ( x, u ) $
where $ f $, $ \partial_t f $, $\partial_u f $ are continuous,
$ \partial_u \tf \geq \b > 0 $ and
$a ( x, u ) $  satisfy  \equ{a1} or \equ{a2}.
Then, for $\e $ small enough,
\equ{equ}-\equ{bc}-\equ{pc} has at least one weak solution $ u \in E $.
If moreover
$ f $, $ \partial_t f$, $ \partial_u f \in C^j (\oO \times \mathbb{R})$,
$ j \geq 1 $, then $ u \in H^{j+1} (\O )\cap C^j_0 (\oO ) $.}

\medskip

In the next subsection  we describe our method of proof.

\subsection{Scheme of the Proof}\label{vlsr}

In order to find critical points of the Lagrangian action functional
$\Psi : E \to \mathbb{R} $ defined in \equ{Psi}
we perform
a variational Lyapunov-Schmidt reduction, decomposing
the space $ E := $ $  H^1( \O )\cap $ $ C^{1/2}_0 (\overline{\O}) $
as
$$
E = V \oplus W
$$
where
$$
V:=N\cap H^1( \O)\qquad \ \mbox{and} \qquad \
W:=N^\bot\cap H^1( \O)\cap C^{1/2}_0 (\overline{\O}) \ .
$$
Setting $ u = v + w $ with $ v \in V$, $ w\in W $ and
denoting by $ \Pi_N $ and $\Pi_{N^\bot}$ the projectors from $L^2 ( \O )
$ onto $ N $ and
$N^\bot$ respectively,
problem \equ{equ}-\equ{bc}-\equ{pc} is equivalent to solve the
{\sl kernel equation} 
\begin{equation}\label{k}
 \Pi_N\, f(v+w,\e)=0
\end{equation}
and the {\sl range equation}
\begin{equation}\label{r}
w = \e \square^{-1}\,\Pi_{N^\bot}\, f(v+w,\e)
\end{equation}
where
$\square^{-1}: N^\bot\to N^\bot$ is the inverse
 of $\square$
and
$f(u,\e)$ denotes the Nemitski operator
associated to $f$, namely
$$
[f(u,\e)](t,x):=f(t,x,u,\e).
$$

\begin{rem}\label{pros}
The usual approach (see [R67]-[DST68]-[T69]-[R71])
is to find, first, by the monotonicity of $ f $, the
unique solution $ v = v( w) $ of the kernel equation \equ{k} and,
next, to solve the range equation \equ{r}.
On the other hand, for nonmonotone forcing terms, one can {\sl not}
in general solve uniquely the kernel equation --recall by
remark \ref{rem:Hp4} that in general unicity
of solutions does not hold.
Therefore we must solve first the range equation
and thereafter the kernel equation.
\end{rem}

We solve,  first,
the range equation by means of a quantitative version of the
Implicit Function Theorem, finding a solution
$ w := w(v, \e ) \in W $ of \equ{r} with
$\| w( v, \e)\|_E = O(\e )$,
see Proposition \ref{pro:range}. Here
no serious difficulties arise since
$ \square^{-1} $ acting on $ W $ is
a compact operator, due to the assumption $T = 2 \pi $, see
\equ{uc}.
\begin{rem}\label{orto}
More in general, $ \square^{-1} $
is compact on the orthogonal complement of ker$( \square)$ whenever
$T$ is a rational multiple of $ 2 \pi $.
On the contrary, if $ T $ is an irrational
multiple of $ 2 \pi $, then
$ \square^{-1} $ is, in general, unbounded
$($a ``small divisor'' problem appears$)$, but
the kernel of
$ \square $ reduces to $ 0 $ $($there is no bifurcation equation$)$.
For existence of periodic
solutions in the case $ T \slash 2 \pi $ is irrational
see [PY89].
\end{rem}

Once the range equation \equ{r} has been solved by $w(v,\e) \in W $
it remains
the infinite dimensional kernel equation
(also called bifurcation equation)
\begin{equation}\label{ke}
\Pi_N f(v+w(v,\e),\e)=0 \ .
\end{equation}

We note (see Lemma \ref{natcon}) that \equ{ke}
is the Euler-Lagrange equation
of the {\sl reduced Lagrangian action functional}
\begin{equation}\label{Phi}
\Phi : V \longrightarrow \mathbb{R}\qquad\qquad
\Phi(v):=\Phi(v,\e):=\Psi(v+w(v,\e),\e) \ .
\end{equation}
$ \Phi $ lacks compactness properties
and, therefore, to find critical points of $ \Phi  $,
we can not rely on critical point theory:

\begin{rem}\label{rem:aut}
Implementing an analogue Lyapunov-Schmidt reduction
in the autonomous case (see [BB03])
it turns out that,
in the corresponding reduced Lagrangian action functional,
a further term proportional to $ \| v \|_{H^1}^2 $ is present.
Therefore it is possible to
apply critical point theory (e.g. the Mountain Pass Theorem)
to find existence and multiplicity of solutions, see [BB04a].
The elliptic term
$ \| v \|_{H^1}^2 $ helps also in proving
regularity results
for the solutions.
\end{rem}

\noindent
We attempt to minimize $ \Phi  $.

\indent
We do not
try to apply the direct methods of the calculus of variations.
Indeed $ \Phi $, even though it could
possess some coercivity property,
will not be convex
(being $f$ non monotone). 
Moreover, without assuming
any growth condition on the nonlinearity $ f $, the functional
$ \Phi $
could neither be well defined on any
$L^p$-space.
\\[1mm]
\indent
Therefore we minimize $\Phi $
in any $ \overline{B_R}
:= $ $ \big\{v\in V, \| v\|_{H^1}\leq R \big\} $, $ \forall R > 0 $,
as in [R67].
By standard compactness arguments
$ \Phi  $ attains minimum at, say, $ \bar v \in \overline{B_R}.$
Since $ \bar v $ could belong to the boundary $ \partial \overline{B_R} $,
$\bv $ could not be a solution of \equ{ke} and
we can only conclude the variational inequality
\begin{equation}\label{minore}
D_v\Phi(\bv)[\f]=\int_\O f(\bv+w(\bv,\e),\e)\f \leq 0
\end{equation}
for any {\sl admissible variation} $ \f\in V $,
i.e. if
$ \bv+\theta\f\in \overline{B_R}$,
$\forall \theta < 0 $ sufficiently small.
\\[1mm]
\indent
The heart of the existence proof of the weak solution $ u $ of
Theorem 1,
Theorem 2 and Theorem 3
is to obtain, choosing suitable admissible variations like in [R67],
the {\sc a-priori estimate}
$\|\bar v\|_{H^1} < R $
for some $ R > 0 $, i.e. to show that
$\bv$ is an {\sl inner} minimum point of $ \Phi $ in $ B_R$.
\\[1mm]
\indent
The strong monotonicity assumption
$ (\partial_u f) (t,x,u) \geq \b > 0 $ would
allow here to get such a-priori estimates by arguments similar to
[R67].
On the contrary, the main difficulty
for proving Theorems
1, 2 and 3
which deal with non-monotone nonlinearities
is to obtain such a priori-estimates for $ \bv $.
\\[1mm]
\indent
The most difficult cases are
the proof of Theorems 1 and 2.
To understand the problem, let consider the particular
nonlinearity  $ f(t,x,u) = u^{2k} + h(t,x) $
of Theorem 1.
The even term $ u^{2k} $ does not give
any contribution into the
variational inequality \equ{minore} at the $ 0^{th} $-order in $\e$,
since the right hand side of \equ{minore} reduces, for $ \e = 0 $, to
$$
\int_\O \Big( \bv^{2k} + h(t,x) \Big) \f =
0, \qquad  \forall \f \in V
$$
by \equ{pari} and  $h \in N^\bot $.

Therefore, for deriving, if ever possible, the required a-priori estimates,
we have to develop the variational inequality \equ{minore}
at higher orders in $ \e $.
We obtain
\begin{equation}\label{varda}
0 \geq  \int_\O 2k \bv^{2k-1} \f \ w(\e, \bv )  +O ( w^2(\e, \bv) )  =
\int_\O \e \ 2k \bv^{2k-1} \f \,
\square^{-1} ( h +\bv^{2k} )  
 + O( \e^2 )
\end{equation}
because
$ w (\bv, \e) = $
$ \e  \square^{-1} ( \bv^{2k} + h ) + o( \e )  $
(recall that $ \bv^{2k} $, $ h \in N^\bot $).

We now sketch
how the $ \e $-order term
in the variational inequality \equ{varda}
allows to prove
an $L^{2k}$-estimate for $ \bv $.
Inserting the admissible variation $ \f := \bv $ in \equ{varda}
we get
\begin{equation}\label{var}
\int_\O H \bv^{2k}  + \bv^{2k} \square^{-1} \bv^{2k}
\leq O (\e)
\end{equation}
where $ H $ is a weak solution of $ \square H = h $ which
verifies $ H ( t, x ) > 0 $ in $ \O $
($ H $ exists by the ``maximum principle'' Proposition \ref{lem:6}).
The crucial fact is that the first term in \equ{var}
satisfies the coercivity inequality
\begin{equation}\label{manwe}
\int_\O H v^{2k} \geq
c(H) \int_\O v^{2k}\,,\qquad \forall \, v\in V 
\end{equation}
for some
constant $ c(H) >0$, see Proposition \ref{pro:fund}.
The second term $ \int_\O \bv^{2k} \square^{-1} \bv^{2k}  $
will be
negligible, $\e$-close to the origin, with respect to
$ \int_\O H v^{2k} $ and \equ{var}, \equ{manwe} will
provide the $L^{2k}$-estimate for $\bv $.

We remark that the inequality  \equ{manwe}
is not trivial because
$H $ vanishes at the boundary
($H(t,0)=$ $H(t,\pi) = 0 $).
Actually,  the proof of
\equ{manwe}
relies on the form $ v(t,x) = $
$\hat v (t+x) -$ $ \hat v (t-x) $
of the functions of $ V $.
\\[1mm]
\indent
Next, we can obtain,
choosing further admissible variations $ \f $
in \equ{varda} and using
inequalities similar to \equ{manwe},
an $ L^\infty $-estimate for $\bv $ and, finally, the
required $H^1 $-estimate,
proving the
existence of a weak solution
$ u \in E $, 
see section \ref{sec:12}.

Moreover, using similar techniques inspired to regularity theory and
further suitable variations, we can also obtain
a-priori estimates for the $L^\infty$-norm of the
higher order derivatives of $\bv $ and for its
$H^j$-Sobolev norms. In this way we can prove
the regularity of the solution $ u  $ 
-fact quite surprising for non-monotone nonlinearities-,
see  subsection \ref{sec:hig}.
\\[1mm]
Theorem 2 is proved developing such ideas
and a careful analysis of the further term $ { \cal R } $.
\\[1mm]
\indent
The proof of Theorem 3 is easier than for
Theorems 1 and 2. Indeed the
additional term $ a(x,u) $ does not contribute into the variational
inequality \equ{minore} at the $ 0^{th} $-order in $ \e $,
because $\int_\O a(x, \bv ) \f \equiv 0 $, $ \forall \f \in V $,
by \equ{pari2}. Therefore the dominant term in the variational
inequality \equ{minore} is provided by the monotone forcing term $ \tf $
and
the required a-priori estimates are obtained
with arguments similar to [R67], see section \ref{sec:thm2}.

\bigskip

{\bf Acknowledgments:}
The authors thank A. Ambrosetti for having suggested the study
of this subject and for useful discussions. They also
thank P. Bolle, G. Mancini, M. Procesi for interesting comments.
Part of this paper was written when the second
author was visiting S.I.S.S.A. in Trieste.
\\[3mm]
{\bf Notations:}
$\O := \mathbb{T} \times (0,\pi) $ where $ \mathbb{T} :=
\mathbb{R} \slash 2 \pi \mathbb{Z} $.
We denote by $ C^j ( \oO ) $ the Banach space of
functions $ u:\oO\longrightarrow \mathbb{R}$
with $j$ derivatives in $ \O $ 
continuous up to the boundary
$ \partial \O $, endowed
with the standard norm $\| \ \|_{C^j}$. $ C^j_0 (\oO ) :=
C^j (\oO ) \cap C_0 (\oO )$ where $C_0 (\oO )$ is
the space of real valued continuous functions satisfying
$u(t,0) =$ $ u(t,\pi)=$ $0$.
Moreover $ H^j (\O) := W^{j,2} (\O) $
are the usual Sobolev spaces with scalar product
$ \langle \cdot , \cdot \rangle_{H^j} $
and norm
$ \| \ \|_{H^j( \O ) }^2 $ .
Here $ C^j ( \mathbb{T} ) $ denotes the Banach space of periodic
functions $ u : \mathbb{T} \to \mathbb{R} $
with  $j$ continuous derivatives.
Finally, $ H^j ( \mathbb{T} ) $ is
the usual Sobolev space of $2 \pi$-periodic functions.

\section{Preliminaries}\label{sec:lemmi}

We first collect some important properties on the D'Alembertian operator
$ \square $.

\begin{df} \label{def1}
Given $ f (t,x) \in L^2 ( \O ) $, a function
$u \in L^2(\O )$ is said to be a weak solution
of $ \square u = f $ in $\O $ satisfying the boundary conditions
$ u(t,0)=$ $ u(t,\pi) =0 $,
iff
$$
\int_\O u \square \f =  \int_\O f \f \qquad \forall \f \in C^2_0 ( \ov{\O} ).
$$
\end{df}

It is easily verified that,
if $u \in C^2 ( \ov{\O} ) $ is a weak solution
of $ \square u = f $ according to definition \ref{def1},
then $u $ is a classical solution and $ u(t,0)=$ $ u(t,\pi) =0 $.

The kernel $ N \subset L^2(\O) $ of the D'Alembertian operator $ \square $,
i.e. the space of weak solutions of the
homogeneous linear equation $ \square v = 0 $
verifying the Dirichlet boundary conditions $ v(t,0)= v(t,\pi) = 0 $,
is the subspace $N$ defined in \equ{KerN}.
$N$ coincides with the closure in $ L^2(\O) $ of the
classical solution of
$ \square v=0 $ verifying Dirichlet boundary conditions
which, as well known, are of the form
$ v (t,x)= $ $\hv (t+x)- $ $\hv(t-x) $, $ \hv\in C^2(\mathbb{T})$.

Using Fourier series we can also characterize $ N $ as
$$
N=\left\{\
 v(t,x) \in L^2 (\O ) \ | \
v(t,x) = \sum_{j \in {\mathbb Z}}
a_j e^{\ii j t } \sin jx \qquad {\rm with} \qquad
\sum_{j\in {\mathbb Z} } |a_j|^2 < \infty   \right\}\, .
$$
The range of $\square $ in $L^2 (\O)$ is
\begin{eqnarray*}
N^\bot &:=& \Big\{ f \in L^2 ( \O ) \ \  \big| \ \
\int_\O f v =0, \  \forall v \in N \Big\} \,,\\
&=&
\Big\{
f(t,x) \in L^2 (\O ) \ \big| \
f(t,x) = \sum_{l \in {\mathbb Z},\, j\geq 1\atop{j\neq |l|}}
f_{lj} e^{\ii l t } \sin jx \quad {\rm with} \quad
\sum_{l \in {\mathbb Z},\, j\geq 1\atop{j\neq |l|}} |f_{lj}|^2 < \infty   \Big\}\, .
\end{eqnarray*}
i.e.
$ \forall  f(t,x)\in N^\bot $ there exists a unique weak solution
$ u = \square^{-1}f \in N^\bot $ of $\square u = f $.

Furthermore
$\square^{-1}$ is a 
bounded operator such that
\begin{equation}\label{W}
\square^{-1}\ :\ N^\bot\longrightarrow
N^\bot \cap H^1(\O)\cap C^{1/2}_0(\oO) \,
\end{equation}
i.e. there exists a suitable constant $\uc\geq 1$ such that
\begin{equation}\label{uc}
\Big\| \square^{-1} f \Big\|_E
\leq \uc\,\| f\|_{L^2} \  \qquad \forall f\in L^2 ( \O )
\end{equation}
where $\| u\|_E := \| u\|_{H^1}+\|u\|_{C^{1/2}}$.
By \equ{uc} and the compact embedding $H^1(\O) \hookrightarrow L^2(\O), $
the operator $\square^{-1}: N^\bot\to N^\bot$ is compact.

These assertions follow easily from the Fourier series representation
(see e.g. [BN78])
$$
f(t,x) :=\sum_{j\geq 1, j\neq |l|} f_{lj} e^{\ii lt } \sin jx \ \
\Rightarrow \ \
u = \square^{-1}f :=
\sum_{j\geq 1, j\neq |l|} \frac{f_{lj}}{-l^2 + j^2}
e^{\ii lt } \sin jx
$$
noting that $ u $ is a weak solution
of \equ{equ} (according to definition \ref{def1})
iff $ u_{lj} = f_{lj} \slash (-l^2 + j^2) $
$ \forall l \in {\bf Z} $, $ j \geq 1 $,
see e.g. [DT68], [H82].
\\[1mm]
\indent
To continue, $ \square^{-1} $ is a bounded operator
also between the spaces
\begin{equation}\label{regdiG}
L^\infty(\O)
\longrightarrow
C^{0,1}(\oO)\,,\qquad
H^k(\O)
\longrightarrow
H^{k+1}(\O)\,,\qquad
C^k(\oO)
\longrightarrow
C^{k+1}(\oO)\,
\end{equation}
as follows by the integral formula
for $ u = \square^{-1} f = $ $ \Pi_{N^\bot} \psi $
where (see e.g. [L69], [BCN80])
$$
\psi (t,x) := - \frac{1}{2} \int_x^\pi  \int_{t+x-\xi}^{t-x+\xi}
f(\xi, \tau) \ d \tau  d \xi + c \frac{(\pi- x)}{\pi}\,,
$$
 $ \Pi_{N^\bot} : L^2 ( \O )\to N^\bot $ is the orthogonal projector
onto $ N^\bot $
and
\begin{equation}\label{clovi}
    c := \frac{1}{2} \int_0^\pi  \int_{t-\xi}^{t+\xi}
f(\xi, \tau) \ d \tau d \xi\equiv
\frac{1}{2} \int_0^\pi  \int_{-\xi}^{\xi}
f(\xi, \tau) \ d \tau d \xi
\equiv \mbox{const}\,
\end{equation}
is a constant independent of $t$, because\footnote{We have that
$ 2 c = $ $ \int_{\mathcal{T}(t)} f = $
$ \lim_{ n \to \infty } \int_{\mathcal{T}(t)} f_n $
where
$ \mathcal{T}(t):= $ $ \{ ( \tau,\xi) \in \O $
s.t. $ t-\xi< \tau<t+\xi$, $0 < \xi < \pi \} $
and
$$
f_n(t,x):=
\sum_{|l|,j\leq n\atop{j\neq |l|}}
f_{l,j} e^{\ii l t}\sin jx\stackrel{L^2}\longrightarrow
\sum_{j\neq |l|} f_{l,j} e^{\ii l t}\sin jx = f(t,x) \in N^\bot \ .
$$
The claim follows since $ \int_{\mathcal{T}(t)} f_n $
is, for any $  n $,  independent on $ t $ :
\begin{eqnarray*}
    \int_{\mathcal{T}(t)} f_n
    &=&
    \int_0^\pi
    \sum_{|l|,j\leq n\atop{j\neq |l|}} f_{l,j}\sin j \xi
    \int_{t-\xi}^{t+\xi} e^{\ii l t} d\tau d\xi\\
    &=&
    \sum_{1\leq j\leq n}f_{0,j}\int_0^\pi 2\xi\sin j\xi\, d\xi
    +
    \sum_{|l|,j\leq n\atop{j\neq |l|,l\neq 0}}
2 \frac{f_{l,j}}{l} e^{\ii l t}  \int_0^\pi \sin j \xi \sin l\xi\, d\xi
    = \sum_{1\leq j\leq n}f_{0,j}\int_0^\pi 2\xi\sin j\xi\, d\xi \ .
\end{eqnarray*}}
$ f\in  N^\bot $.

We also have,
since $ \partial^j_t H  $ is a weak solution of
$ \square (\partial^j_t H ) = \partial^j_t h $ and \equ{W} applies,
\begin{equation}\label{dertj}
h \in H^j ( \O ) \qquad \Longrightarrow \qquad \partial^j_t H \in
C^{1\slash2}_0 ( \oO )\ .
\end{equation}
Finally, the projector $ \Pi_N : L^2 (\O) \to N $ can be written
as $\Pi_N u = p(t+x ) - p(t-x) $ where
$$
p(y) := \frac{1}{2\pi} \int_0^\pi \Big[ u (y-s,s) -
u (y +s,s) \Big] ds \
$$
and therefore, since $ u \in C^j ( \oO) \Rightarrow p \in C^j ( \mathbb{T}) $
and $ u \in H^j ( \O) \Rightarrow p \in H^j ( \mathbb{T}) $,
\begin{equation}\label{Lovip}
\Pi_N, \ \Pi_{N\bot}
: C^j (\oO ) \to  C^j (\oO) \qquad \mbox{are bounded}
\end{equation}
\begin{equation}\label{proHj}
\Pi_N, \ \Pi_{N\bot}
: H^j ( \O ) \to  H^j ( \O) \qquad \mbox{are bounded}.
\end{equation}

\medskip\medskip

\noindent
{\bf Kernel properties and technical Lemmata}

\medskip\medskip

Let define, for  $0\leq \a<1/2 $,
\begin{equation}\label{Oa}
\O_\a:=\mathbb{T}\times(\a\pi,\pi-\a\pi)\subseteq\O \ .
\end{equation}
\begin{lem}\label{lem:0}
Let $a\in L^1(\O)$. 
\begin{equation}\label{fubini1}
\int_{\O_\a}a(t,x)\,dt\,dx =
\frac{1}{2}\intt \int_{-2\pi+s_++2\a\pi}^{s_+-2\a\pi}
a\left( \frac{s_+ + s_-}{2}, \frac{s_+ - s_-}{2}\right)
\, ds_-\, ds_+\ .
\end{equation}
In particular for $p,q\in L^1(\mathbb{T})$
\begin{equation}\label{fubini2}
    \intOa p(t+x) q(t-x)\,dt\,dx=\frac{1}{2}\intt \!\!p(s)ds \intt\!\! q(s)ds-
    \frac12 \int_{-2\a\pi}^{2\a\pi} \intt\!\! p(y)q(z+y)dydz
\end{equation}
and
\begin{equation}\label{fubini6}
    \int_{\O_\a} p(t+x)\,dt\,dx =\int_{\O_\a} p(t-x)\,dt\,dx = \pi (1-2\a)\intt p(s)ds\ .
\end{equation}
Moreover,
given $f,g:\mathbb{R}\rightarrow \mathbb{R}$ continuous,
\begin{equation}\label{fubini7}
\int_{\O_\a} f(p(t+x)) g(p(t-x))\,dt\,dx =
\int_{\O_\a} f(p(t-x)) g(p(t+x))\,dt\,dx\ .
\end{equation}
\end{lem}

\textsc{Proof}: In the Appendix. \qed

\medskip\medskip

\begin{lem}\label{iso}
For any $ v = $ $\hv(t+x)-\hv(t-x) := $ $ v_+ - v_-  \in N $
\begin{equation}\label{isometria}
    \| v \|^2_{L^2(\O)}=2\pi \| \hv \|^2_{L^2(\mathbb{T})}=
     2\pi\int_0^{2\pi}\hv^2 \ .
\end{equation}
Moreover
\begin{equation}\label{isometriaH}
  \| v_t \|^2_{L^2(\O)}=\| v_x \|^2_{L^2(\O)}=
2\pi \| \hv' \|^2_{L^2(\mathbb{T})}\ \qquad \forall v \in N
\cap H^1 (\O )
\end{equation}
\begin{equation}\label{isometriaL}
\| \hv \|_{ L^{\infty} (\mathbb{T} ) } \leq \| v \|_{L^\infty(\O)}
\leq  2\| \hv \|_{L^\infty (\mathbb{T})} \ \qquad
\forall v \in N \cap L^\infty (\O)
\end{equation}
\begin{equation}\label{embed}
\| v \|_{L^\infty(\O)} \leq \| v \|_{H^1(\O)}\  \
\qquad \qquad \forall v \in N \cap H^1 ( \O ) \ .
\end{equation}

\end{lem}

\textsc{Proof}: In the Appendix. \qed

\medskip\medskip

\begin{lem}\label{lem:2}
Let $\f_1,\dots,\f_{2k+1}\in N $ 
and $\f_1\cdot\ldots\cdot \f_{2k+1}\in L^1(\O).$
Then 
\begin{equation}\label{fubini4}
   \intOa \f_1\cdot\ldots\cdot \f_{2k+1}=0\ .
\end{equation}
In particular $\f_1\cdot\ldots\cdot \f_{2k}\in N^\bot.$

Moreover, if
$ a: \O \longrightarrow\mathbb{R} $ satisfies
($i$) $ a(x,u) = a ( \pi - x, u) $ and $ a ( x, u ) = a(x,-u) $
or ($ii$) $ a(x,u) = - a( \pi - x, u) $ and $ a ( x, u ) = - a(x,-u) $,
then
\begin{equation}\label{pari}
\intOa a(x, v)\,\f =0 \qquad \forall\, v\in N\cap L^\infty\ , \ \f\in N
\end{equation}
and
\begin{equation}\label{pari2}
\intOa (\partial_u a) (x, v)\, \f_1 \f_2  =0
\qquad \forall\, v\in N\cap L^\infty\ , \ \f_1, \f_2 \in N\ .
\end{equation}
\end{lem}

\textsc{Proof}: In the Appendix. \qed

\medskip\medskip

\begin{lem}\label{lem:ineq}
The following inequalities hold:
\begin{equation}\label{dis:ovvia}
(a-b)^{2k}\leq 2^{2k-1}(a^{2k}+b^{2k}) \qquad \ \  \
\forall\ a,b\in \mathbb{R}\ ,
\end{equation}
\begin{equation}\label{zialalletta}
(a-b)^{2k}\geq a^{2k}+b^{2k}-2k(a^{2k-1}b+ab^{2k-1}) \qquad
\forall\  a,b\in \mathbb{R},\ \
k\in \mathbb{N},\, k\geq 2
\end{equation}
\begin{equation}\label{zialalletta3}
(a+b)^{2k-1}-a^{2k-1}\geq 4^{1-k}\,b^{2k-1}
\qquad\qquad\forall \ a\in\mathbb{R}\ ,\ b>0\ ,\
 k\in\mathbb{N^+}\ .
\end{equation}
\end{lem}

\textsc{Proof}: In the Appendix. \qed

\begin{lem}\label{lem:3}
Let $v \in N \cap L^{2k} ( \O )$ and $k\in \mathbb{N}^{+}$. Then
\begin{equation}\label{fubini5}
\int_\O v^{2k}\leq \pi 4^{k}\intt \hv^{2k} \ .
\end{equation}

Moreover, for $k=1$
$$
\int_{\O_{\a}} v^{2} \geq  2\pi \Big[1- 4\a\Big]\intt \hv^{2}
\geq \pi \intt \hv^{2}
$$
if $ 0 \leq \a\leq 1 \slash 8 $. For $ k \geq 2, $
$$
\int_{\O_{\a}} v^{2k} \geq  2\pi \Big[1- 2(1+2k)\a\Big]\intt \hv^{2k}
\geq \pi \intt \hv^{2k}
$$
if $0 \leq \a\leq 1 \slash 4(1+2k)$.
\end{lem}

\textsc{Proof}: In the Appendix. \qed

\medskip\medskip

\noindent
{\bf Generalities about the difference quotients}

\medskip

\noindent
For $f\in L^2 (\O)$ we define the difference quotient of size
$ h \in\mathbb{R}\setminus\{ 0 \}$
$$
(D_h f) (t,x):=\frac{f(t+h,x)-f(t,x)}{h}\,
$$
and the $h$-translation
$$
(T_h f) (t, x):=f(t+h, x )
$$
with respect to time.

The following Lemma
collects some elementary properties of the difference quotient.

\begin{lem}\label{lem:quoz}
Let  $f,g\in L^2 (\O)$, $h\in\mathbb{R}\setminus\{ 0 \}$.
The following holds
\begin{itemize}
\item[(i)] Leibniz rule:
\begin{eqnarray}\label{Leibniz}
 D_h (fg) &=& (D_h f)g+T_h f D_h g\ ,\\
 D_h f^m  &=& (D_h f) \sum_{j=0}^{m-1} f^{m-j-1} T_h f^j\nonumber\\
          &=& m(D_h f)f^{m-1}  + (D_h f)
                 \sum_{j=0}^{m-1} f^{m-j-1} (T_h f^j-f^j)\label{Leibniz2}\ ,\\
\int_\O D_h (fg) &=&\int_\O (D_h f)g+ f (D_{-h} g)\ ;\label{Leibnizint}
\end{eqnarray}
\item[(ii)] integration by parts:
\begin{equation}\label{perparti}
    \int_\O f (D_{-h}g)=-\int_\O (D_h f) g\ ;
\end{equation}
\item[(iii)] weak derivative: If there exists a constant $ C > 0 $
such that $\forall h $ small
\begin{equation}\label{weakd}
\| D_{h} f \|_{L^2 } \leq C \quad \Rightarrow \quad
{\rm then \ } f { \ \rm has \ a \ weak \ derivative \ } f_t  \ {\rm and \ } \
\| f_t \|_{L^2} \leq C
\end{equation}
\end{itemize}
Moreover, if $ f $ has a weak derivative $f_t  \in L^2 (\O) $, then
\begin{itemize}
\item[(iv)] estimate on the difference quotient:
\begin{equation}\label{sti:diffquo}
\| D_h f  \|_{L^2(\O)} \leq
\|  f_t  \|_{L^2(\O)}\ ;
\end{equation}
\item[(v)] convergence:
\begin{equation}\label{conL2quoz}
D_h f \stackrel{L^2}\longrightarrow f_t\  \qquad\mbox{as}
\quad h\longrightarrow 0\ .
\end{equation}
\end{itemize}
\end{lem}

\noindent
\textsc{Proof}. In the Appendix. \qed

\section{The Lyapunov-Schmidt decomposition}\label{sec:L-S}

\subsection{The Range equation}\label{sec:range}

We first solve the range equation \equ{r}
applying the following quantitative version
of the Implicit Function Theorem, whose standard proof is omitted.

\begin{pro}\label{TFI}
Let $X,Y,Z$ be Banach spaces and $x_0\in X$, $y_0\in Y$.
Fix $r, \r>0$ and  define $X_r:=\{ x\in X\ \mbox{s.t.}\  \| x-x_0\|_X< r  \}$
and $Y_\r:=\{ y\in Y
\ \mbox{s.t.}\  \| y-y_0\|_Y< \r  \}$.
Let ${\cal F} \in \mathcal{C}^1(U,Z)$
where $\overline{Y_\r}\times \overline{X_r}
\subset U\subset Y\times X$
is an open set.
Suppose that
\begin{equation}\label{TFI1}
    {\cal F}(y_0,x_0)=0
\end{equation}
and that $D_y {\cal F}(y_0,x_0)\in \mathcal{L}(Y,Z)$ is invertible.
Let $T:=\Big( D_y {\cal F}(y_0,x_0) \Big)^{-1}$ and $\|T\| :=
\|T\|_{\mathcal{L}(Z,Y)} $ be its norm.
If
\begin{eqnarray}
   \sup_{\overline{X_r}}\| {\cal F}(y_0,x) \|_Z &\leq& \frac{\r}{2\| T\|}\label{TFI2}\,, \\
   \sup_{\overline{Y_\r}\times \overline{X_r}}
   \Big\| {\rm Id}_Y-T D_y {\cal F}(y,x)
\Big\|_{\mathcal{L}(Y,Y)} &\leq& \frac12\,,\label{TFI3}
\end{eqnarray}
then
there exists $y\in \mathcal{C}^1(X_r,\overline{Y_\r})$
such that ${\cal F}(y(x),x)\equiv 0.$

\end{pro}

Appling Proposition \ref{TFI} to the range equation \equ{r} we derive:

\begin{pro}\label{pro:range}
Let $f=f(t,x,u,\e)$, $f_u :=\partial_u f $ and
$\e \partial_\e f $ be continuous on $\oO\times\mathbb{R}\times[-1,1]$.
Then $ \forall\, R>0 $ there exists a unique function
\begin{equation}\label{wC0}
w=w(v,\e)\in \mathcal{C}^1\bigg(
\{ \|v\|_{L^\infty}< 2 R  \} \times
\{ |\e|<\e_0(R)  \}\, ,\, \Big\{ \| w\|_E \leq C_0(R)|\e | \Big\}
 \bigg)\,
\end{equation}
solving the range equation \equ{r},
where $\e_0(R):=1/ 2C_0(R)$ and\footnote{$\uc$ is defined in \equ{uc}.}
\begin{equation}\label{C0}
    C_0(R):=1+\sqrt{2}\pi \uc
\max_{\oO\times\{|u|\leq 3 R + 1 \}\times\{ |\e|\leq 1 \}}
    \Big[ |f(t,x,u,\e)|+ |f_u(t,x,u,\e)|  \Big] \ .
\end{equation}

\noindent
Moreover, the following continuity property holds
\begin{equation}\label{vnv*}
v_n\,,\ \bv\in V\,,\ \  \|v_n\|_{L^\infty}\, ,\ \|\bv\|_{L^\infty}\leq R\, ,\ \
v_n\stackrel{L^2}\longrightarrow \bv\ \ \ \Longrightarrow\ \ \
w(v_n)\stackrel{E}\longrightarrow w(\bv)\,.
\end{equation}

\noindent
\end{pro}
\textsc{Proof}:
Let
$X :=V\times \mathbb{R},$ $Y = Z := W$ (namely $x:=(v,\e)$ and $y:=w$)
and $\| x\|_X:= $ $ \| v\|_{L^\infty}+\frac{R}{\e_0(R)}|\e|$.
Let also $x_0:=(0,0),$ $y_0:=0$, $r:=3R,$ $\r:=1$,
${\cal F}(y,x):= $ ${\cal F}(w,v,\e):= $
$ w-\e \square^{-1}\Pi_{N^\bot}f(v+w,\e)$
and $U:= $ $ W\times V\times (-1,1)$. Note that ${\cal F}( \cdot , \cdot )
\in C^1 $ since  the Nemitski operator
$ \e f \in $ $ C^1 (E \times (-1,1), L^2 (\O)) $ and \equ{W} holds.
Moreover \equ{TFI1} holds and $T={\rm Id}_W$ (hence $\| T\|=1$).

If
$\|v\|_{L^\infty} \leq 3R $
and $\| w\|_{L^\infty}\leq \| w\|_E \leq 1$, then

\begin{equation}\label{R+1}
|v(t,x)+w(t,x)|\leq 3R+1\,,\qquad\forall\, (t,x)\in\oO\,.
\end{equation}
\noindent
Using \equ{uc}, Bessel inequality
$\|\Pi_{N^\bot}f\|_{L^2}\leq \| f\|_{L^2} $,
$ \| f \|_{L^2(\O)}\leq \sqrt{2}\pi\| f \|_{L^\infty(\O)} $ and
estimate \equ{R+1}, we obtain, $ \forall |\e | \leq \e_0 (R) $
and $ \|v\|_{L^\infty} \leq 3R $,
\begin{equation}
    \| {\cal F}(0,v,\e)\|_E \leq  |\e| \ \uc \| f(v,\e) \|_{L^2}
    \leq |\e | \ \sqrt{2}\pi \uc \| f(v,\e) \|_{L^\infty}\leq C_0(R) |\e |
\leq \frac12
\label{ettore}
\end{equation}
where $C_0 (R) $ is defined in \equ{C0}.
Hence
\equ{TFI2} follows from \equ{ettore}.

\noindent
Since
$$
D_w {\cal F}(w,v,\e)[\wtilde w]=\wtilde w -\e \square^{-1}\Pi_{N^\bot}
\Big( f_u (v+w,\e)  \wtilde w \Big)\qquad \quad \forall\, \wtilde w\in W \, ,
$$
we deduce, arguing as before, $\forall \|v\|_{L^\infty} \leq 3R $,
$\forall \| w \|_E \leq 1 $, $ \forall |\e | \leq \e_0  $,
\begin{equation}
   \sup_{\|\wtilde w\|_E=1} \Big\| \wtilde w -
D_w {\cal F}(w,v,\e)[\wtilde w] \Big\|_E \leq
|\e| \ \sqrt{2}\pi \uc \| f_u(v+w,\e) \|_{L^\infty}
\leq
\e_0(R)\, C_0(R) = \frac12
\end{equation}
and \equ{TFI3} follows.
Now we can apply Proposition \ref{TFI}
finding a function
$ w = w ( v, \e ) \in $ $ C^1 ( \{ \|x\|_X< r \}, \ov{W}_1  ) $
satisfying the range equation \equ{r}.
Finally, note that
$$
\{ \|v\|_{L^\infty}< 2 R  \} \times
\{ |\e|<\e_0(R)  \}\subset  \{ \|x\|_X< r=3R  \}
$$
and, arguing as above,
$$
\Big|\Big| w (v, \e)  \Big|\Big|_E = \Big|\Big|
 \e  \ \square^{-1} \Pi_{N^\bot} f (v + w(v,\e);
\e) \Big|\Big|_E  \leq 
|\e | \sqrt{2}\pi \uc \| f(v + w; \e)
\|_{L^\infty}\leq |\e | C_0(R) \ ,
$$
whence \equ{wC0} follows.

\noindent
We now prove \equ{vnv*}.
Let $w_n:=w(v_n,\e) $ and $\bar w:=w(\bv,\e) $ denote, for brevity,
the solutions
of $ w_n = \e \square^{-1} \Pi_{N^\bot} f(v_n+w_n)$  and
$ \bar w = \e \square^{-1} \Pi_{N^\bot} f( \bar v + \bar w) $ .
By \equ{uc}, using that  $\partial_u f $ is continuous in
$\oO\times\mathbb{R}\times[-1,1]$,
recalling the definition of $C_0 (R)$ in \equ{C0} and
 that $ \e_0 (R) C_0 (R) = 1 \slash 2 $
\begin{eqnarray*}
   \| w_n-\bar w\|_E &\leq&
   | \e | \  \uc \Big\| f(v_n+w_n,\e)-f(\bv+\bar w,\e) \Big\|_{L^2} \\
   &\leq& \e_0(R)\,\uc \ \| (v_n-\bv)+(w_n-\bar w) \|_{L^2}
\max_{\oO\times\{|u|\leq 3 R
+ 1 \} \times\{ |\e|\leq 1 \}} |f_u (t,x,u,\e)|
 \\
&\leq& \e_0(R)\, C_0 (R)  \| v_n-\bv \|_{L^2}
+ \e_0(R)\, C_0 (R)  \| w_n -\bar w\|_{L^\infty}\\
& = & \frac12 \| v_n-\bv \|_{L^2}
+ \frac12 \| w_n -\bar w\|_E\,,
\end{eqnarray*}
whence
$ \| w_n-\bar w\|_E \leq $ $ \| v_n-\bv \|_{L^2} $
and \equ{vnv*} follows.
\qed

\subsection{The Kernel equation}\label{sec:kernel}

Once the range equation \equ{r} has been solved by $w(v,\e) \in W $
there remains the infinite dimensional kernel equation \equ{ke}.

Since $ V $ is dense in $ N $ with the $L^2$-norm,
equation \equ{ke} is equivalent to
\begin{equation}\label{ker}
\int_\O f(v+w(v,\e),\e)\f =0\, \qquad \forall\, \f\in V\
\end{equation}
which is the Euler-Lagrange equation of the {\sl reduced
Lagrangian action functional}
$\Phi :$ $ V \longrightarrow \mathbb{R} $,
$ \Phi(v):= $
$\Phi(v,\e):= $ $ \Psi(v+w(v,\e),\e) $, defined in \equ{Phi}.
Indeed:

\begin{lem}\label{natcon}
$ \Phi \in $ $ \mathcal{C}^1(\{ \|v\|_{H^1}< 2 R  \} , \mathbb{R})$ and a
critical point $ \bv $ of $\Phi $ is a weak solution of the kernel equation
\equ{ke}.  Moreover $\Phi $ can be written as
\begin{equation}\label{Phi*}
\Phi(v) =
\e\int_\O \Big[ F(v+w(v);\e)-\frac12 f(v+w(v);\e )w(v)\Big] \ dt dx \ .
\end{equation}
\end{lem}

\noindent
\textsc{Proof}.
Since $ \Psi ( \cdot, \e ) \in C^1 (E , \mathbb{R} )$ and,
by Proposition \ref{pro:range}, $w (\cdot, \e ) \in $
$C^1 ( \{ \|v\|_{H^1} < 2 R  \} , \mathbb{R} )$ (note that
$\{ \|v\|_{H^1} < 2 R \} \subset$ $\{ \|v\|_{L^\infty}
< 2 R \}$ by \equ{embed}),
then
$\Phi \in$   $ \mathcal{C}^1(
\{ \|v\|_{H^1}< 2 R  \} , \mathbb{R})$ and
\begin{equation}\label{dfrid1}
D\Phi(v)[\f]=D\Psi(v+w(v)) \Big[ \f+dw(v)[\f] \Big] \quad
\qquad \forall\, \f \in V  \ .
\end{equation}
We claim that,
since $ w = w(v) \in E $ is a weak
solution of the range equation \equ{r}
and $ \tw := dw(v)[\f] \in W $, then
\begin{equation}\label{diffn}
D\Psi (v +w(v) ) \Big[ dw(v)[\f] \Big]=0 \ .
\end{equation}
Indeed, since
$ v_t, v_x \in N $ and 
$ \tw_t$, $ \tw_x \in N^\bot $,
\begin{eqnarray}
 D \Psi(v+w)[\tw] &=&
\int_\O (v+w)_t \tw_t-(v+w)_x\tw_x+\e f(v+w,\e)\tw \label{patroclo}  \\
   & = &
\int_\O w_t \tw_t-w_x\tw_x+\e \Pi_{N^\bot}f(v+w,\e)\tw =  0  \nonumber
\end{eqnarray}
because $ w \in E $ is a weak
solution of the range equation
$\square w=\e \Pi_{N^\bot}f(v+w,\e)$ and 
$w(t,0) = w(t , \pi) =0$.

By (\ref{dfrid1}), \equ{diffn} and
since $ w_t, $ $w_x\in N^\bot $ and $\f_t, $ $ \f_x \in N$
\begin{eqnarray}
  D\Phi(v)[\f] &=& D\Psi(v+w)[\f]
   =   \int_\O (v+w)_t \f_t-(v+w)_x\f_x+\e f(v+w,\e)\f \nonumber\\
   &=& \int_\O v_t \f_t-v_x\f_x+\e f(v+w,\e)\f \label{menelao} \\
   &=& \e\int_\O  f(v+w,\e)\f\, = \e\int_\O  \Pi_{N^\bot} f(v+w,\e)\f\,
\nonumber
\end{eqnarray}
where
in \equ{menelao} we used $ \int_\O v_t \f_t-v_x\f_x=0$ since
$ v, \varphi \in V $.

Finally we prove \equ{Phi*} as in [BB03]. Since
$v_t, v_x \in N $, $w_t, w_x \in N^\bot $ and \equ{isometriaH}
\begin{eqnarray*}
\Phi (v ) & = & \int_\O \frac{(v+w(v))_t^2}{2} -
\frac{(v+w(v))_x^2}{2} + \e F(v + w(v); \e) \\
& = &  \int_\O \frac{(w(v))_t^2}{2} -
\frac{(w(v))_x^2}{2} + \e F(v + w(v); \e)
\end{eqnarray*}
and since
$ \int_\O (w(v))_t^2 -
(w(v))_x^2  = - \int_\O \e f(v + w(v);\e) w(v) $
we deduce \equ{Phi*}. \qed

\medskip\medskip

The next Lemma proves a $L^2$-continuity property for $\Phi $.
\begin{lem}\label{conti}
Let $R>0$ and $|\e|\leq \e_0(R)$ (where
$ \e_0(R) < 1 $ is defined in Proposition
\ref{pro:range}). Then
\begin{equation}\label{continuitaL2}
   v_n,\, \bv\in V \ ,\quad \|v_n\|_{H^1}, \|\bv\|_{H^1}\leq R\ ,\quad
   v_n \stackrel{L^2}\longrightarrow \bv\qquad \Longrightarrow\qquad
   \Phi(v_n)\longrightarrow \Phi( \bv)\ .
\end{equation}
\end{lem}
\textsc{Proof}:
Setting $ w_n:=w(v_n,\e)$ and $\bar w:=w(\bv,\e)$,
we have
\begin{eqnarray*}
\left|
\int_\O F(v_n+w_n)-F(\bv +\bar w)
\right|
   &\leq&
\max_{\oO\times\{|u|\leq R+1 \}\times\{ |\e|\leq 1 \}}
     |f(t,x,u,\e)| \int_\O
     |v_n-\bv + w_n -\bar w|
    \\
   &\leq&
C_0(R)\Big( \|v_n-\bv\|_{L^1}+\|w_n-\bar w\|_{L^1} \Big)
\longrightarrow 0
\end{eqnarray*}
as $n\to\infty$,
 by \equ{vnv*} and the fact that $\| v_n -\bv\|_{L^2}\to 0$.
An analogous estimate holds for the second term in the integral in
\equ{Phi*} proving the Lemma.
\qed

\medskip\medskip

By standard compactness argument the functional
$ \Phi $ attains minimum (resp. maximum) in $ \overline{B_{R}}
:= $ $ \big\{v\in V, \| v\|_{H^1}\leq {R} \big\}$,
$ \forall R > 0 $. Indeed,
let $ v_n \in \overline{B_{R}}$ be a minimazing (resp. maximizing) sequence
$ \Phi(v_n) $ $\rightarrow$ $\inf_{\overline{B_{R}}}\Phi$.
Since $\{ v_n\}_{n\in\mathbb{N}}$ is bounded in $ N \cap H^1$,
up to a subsequence
$ v_n\stackrel{H^1}\rightharpoonup\bar v$
for some $\bv=\bv({R},\e) \in \overline{B_{R}} $.
By the Rellich Theorem
$v_n\stackrel{L^2}\rightarrow\bar v $ and therefore,
by Lemma \ref{conti},
$ \bar v $ is a minimum (resp. maximum) point of $\Phi $ restricted to $\overline{B_{R}}.$
\\[1mm]
\indent
Since $ \bar v $ could belong to the boundary $\partial \overline{B_{R}}$
we only have the variational inequality \equ{minore}
for any {\sl admissible variation} $ \f\in V $,
namely for any $ \f\in V $ such that $\bv+\theta\f\in \overline{B_{R}}$,
$\forall \theta < 0 $ sufficiently small.
As proved by Rabinowitz [R67], a sufficient condition for $ \f \in V $ to be an
admissible variation is the positivity of the scalar product
\begin{equation}\label{scalar}
    \langle \bv,\f  \rangle_{H^1} >0\ .
\end{equation}

\noindent
The heart of the existence proof of Theorems
1, 2 and 3
is to obtain, choosing suitable admissible variations,
the \textsl{a-priori estimate}
$\|\bar v\|_{H^1} < {R} $ for some $ R > 0 $.
i.e. to show that
$\bv$ is an \textsl{inner} minimum (resp. maximum) point of
$ \Phi $ in $ B_{R}$.

\section{Proof of Theorems 1 and 2}\label{sec:12}

The main difficulty for proving Theorems 1 and 2
is
to obtain the fore mentioned a priori-estimate for $\bv $.

\subsection{Proof of Theorem 2}

We look for small amplitude solutions of
(\ref{equ}) with forcing term $ f(t,x,u)$ $
= $ $ g(t,x,u) +$ $ h(t,x) $ where 
$ g (t, x, u ) =$ $ \b(x) u^{2k} +\mathcal{R}(t,x,u)$,
$ h (t,x) \in N^\bot $
and $\mathcal{R}(t,x,u) $ satisfies \equ{h}.
\\[1mm]
Perform the change of variables $ u =  \e ( H+\wtilde u ) $
and set $ \wtilde \e:=\e^{2k} $
\begin{eqnarray*}
\square \wtilde u =
g\big(t, x, \e(H+\wtilde u)\big)
& = & \e^{2k} \Big[
\b(x) (H+\wtilde u)^{2k}+\e^{-2k}\mathcal{R}\big(t, x, \e(H+\wtilde u) \big)
\Big] \\
& = &
\wtilde\e\Big[
\b(x) (H+\wtilde u)^{2k}+
\wtilde\e^{-1}\mathcal{R}\big(t, x, \wtilde\e^{\frac{1}{2k}}(H+\wtilde u) \big)
\Big] \ .
\end{eqnarray*}
Recalling $ \wtilde u \to u $, $ \wtilde\e \to \e $,
we look for solutions of the problem
\begin{equation}\label{equbis}
\begin{cases}
\square u = \e f ( t, x, u; \e ) \cr
u ( t, 0 ) = u ( t, \pi ) = 0 \cr
u(t+2\pi, x ) = u(t, x)
\end{cases}
\end{equation}
where the nonlinear forcing term is
\begin{equation}\label{fnuova}
f(t,x,u;\e) :=
\b(x)(H+ u)^{2k}+\mathcal{R}^*(t,x,u;\e)
\end{equation}
and
\begin{equation}\label{R**}
    \mathcal{R}^*(t,x,u;\e):=
\e^{-1}\mathcal{R}\Big( t,x,\e^{\frac{1}{2k}}\big(H(t,x)+ u\big) \Big)\ .
\end{equation}
Moreover, eventually substituting $\e\to -\e$ and $\b\to -\b$
we can always suppose
\begin{equation}\label{bH}
    \b(x)H(t,x)>0\,,\qquad\qquad
    \forall\, (t,x)\in \O\,.
\end{equation}
By \equ{h},
$\mathcal{R}^*,$
$\partial_u\mathcal{R}^*$,
$ \e \partial_\e \mathcal{R}^* $
are continuous in $\oO \times \mathbb{R} \times [-1,1]$
(recall that $ H \in E$) and
\begin{equation}\label{R**u}
\forall\, R_0 > 0\,,\qquad\qquad
\|\mathcal{R}^*(\cdot;\e)\|_{C(\oO\times\{|u|\leq R_0\})} \,,
\|\partial_u \mathcal{R}^*(\cdot;\e)\|_{C(\oO\times\{|u|\leq R_0\})}
    \stackrel{\e\to 0}\longrightarrow 0 \ .
\end{equation}
Moreover, since $H\in H^1 (\O )$,
$\partial_t \mathcal{R}^*=$
$\e^{-1}\partial_t \mathcal{R}+$
$\e^{\frac{1}{2k}-1}\partial_u \mathcal{R}\  H_t$
and \equ{h}, then
$$
\Big\| \partial_t \mathcal{R}^*(\cdot,u(\cdot);\e) \Big\|_{L^2(\O)}
 \leq
 C^*(\| u\|_{L^\infty(\O)})\qquad
 \forall\ u\in L^\infty(\O)\,,
$$
for a suitable increasing function $C^*(\cdot)$.
\\[1mm]
\indent
In order to find solutions of problem
\equ{equbis} we perform the Lyapunov-Schmidt
reduction of the previous section.
We fix $R>0$ to be chosen later (large enough!).
Since $f, \partial_u f, \e\partial_\e f$ are continuous on
$\oO\times\mathbb{R}\times [-1,1]$,
using Proposition \ref{pro:range}  we solve the range equation
\equ{r} finding $w=w(v,\e)$ for $\| v\|_{L^\infty}<2R$
and $|\e|<\e_0:=$ $\e_0(R)$.
Now we look for minimum or maximum points of
  the corresponding
reduced action functional $\Phi $ in $ \ov{B_{R}} $
according to whether $\e>0$ or $\e<0.$
Since $\Phi$
 attains minimum or maximum  at some point
$ \bv := $ $\bv(\e):= $ $\bv(t,x;\e)$
in $ \ov{B_{R}} $,
to conclude the existence
proof of Theorem 2-(i) we need to show that
$ \bv $ is an interior point in $ B_{R}$, i.e.
$\| \bv \|_{H^1} < R $, for a suitable choice of $R$ large enough.
Let $ \bw:= $ $\bw(t,x;\e):=$ $w(\bv(\e) ,\e)(t,x)\in E$
and $\bu:=$ $\bu(t,x;\e):=$ $\bv+\bw\in E$.
By \equ{wC0} and the definitions of $C_0(\cdot),$ $\e_0(\cdot)$
given in Proposition \ref{pro:range},
 we have
\begin{equation}\label{wC0*}
    \| \bw \|_E\leq C_0(R)|\e |\leq \frac12\,,\qquad\qquad
    \| \bu\|_{C(\oO)}\leq R+\frac12\,,
\end{equation}
since, by \equ{embed}, $\|\bv\|_{C(\oO)}\leq \|\bv\|_{H^1(\O)}\leq R.$
Let
\begin{equation}\label{R*}
     \mathcal{R}_*(t,x;\e):=
     \mathcal{R}^*(t,x,\bu(t,x;\e);\e)
     \,.
\end{equation}
We have $\mathcal{R}_*\in C(\oO)$ and,
choosing $R_0 := R + 1/2 $ in \equ{R**u},
\begin{equation}\label{R*u}
\|\mathcal{R}_*(\cdot;\e)\|_{C(\oO)} \,
    \stackrel{\e\to 0}\longrightarrow 0\,.
\end{equation}
Moreover, since
$\partial_t \mathcal{R}_*=$
$\partial_t \mathcal{R}^*(\bu)+$
$\partial_u \mathcal{R}^*(\bu)\  (\bv_t+\bw_t)$,
we have $\partial_t \mathcal{R}_*$ $\in L^2(\O)$ with
\begin{equation}\label{R*t}
\Big\| \partial_t \mathcal{R}_*(\cdot;\e) \Big\|_{L^2(\O)}
\leq C_*\big(\| \bv\|_{L^\infty(\O)}\big) +o(1)\| \bv_t\|_{L^2(\O)}\,,
\end{equation}
for a suitable increasing function $C_*(\cdot)$.
By \equ{fnuova}, \equ{R**} and \equ{R*}, the variational
inequality \equ{minore} yields, for any admissible variation $ \f \in V $,
$$
0\geq \e \int_\O \b(x) \big(H+ \bv+\bw\big)^{2k}\f + \mathcal{R}_*\f \qquad\qquad
\forall\,\e>0
$$
if $\bv$ is a minimum point,
respectively
$$
0\leq \e \int_\O \b(x) \big(H+ \bv+\bw\big)^{2k}\f + \mathcal{R}_*\f \qquad\qquad
\forall\,\e<0
$$
if $\bv$ is a maximum point.
However, in both cases we get, dividing by $\e$,
\begin{equation}\label{minore0}
\int_\O \b(x) \big(H+ \bv+\bw\big)^{2k}\f \leq
\int_\O - \mathcal{R}_*\f\,.
\end{equation}

The required a-priori estimate for the $ H^1 $-norm of $ \bv $
will be proved in several steps
inserting into the variational inequality
\equ{minore0} suitable admissible variations.
We shall derive, first,
an $ L^{2k} $-estimate for $ \bv $
(it is needed at least when $ k \geq 2 $), see \equ{sti:L2k_fin},
next, an $ L^\infty $-estimate,
see \equ{sti:Li_fin}, and, finally,
the $H^1$-estimate, see \equ{sti:H1_fin1}.
\\[1mm]
\indent
The following key estimate will be heavily exploited.

\begin{pro}\label{pro:fund}
Let $ k \in \mathbb{N}^+ $ and
$ B  \in C (\oO) $ with $ B \geq  0 $ in $ \O$. Define
\begin{equation}\label{sti:costfund}
c_k(B):=\frac{1}{ 4^k}\min_{\bar\O_{\a_k}} B  \
\qquad \mbox{where} \quad  \a_1:=\frac18\ ,\ \ \
\a_k:=\frac{1}{4(1+2k)}\  \ \  \mbox{for}\  \ k\geq 2\ .
\end{equation}
Then
$\forall\ v\in N\cap L^{2k}(\O)$
\begin{equation}\label{sti:fund}
\int_\O B v^{2k} \geq c_k(B) \int_\O v^{2k}\,.
\end{equation}
\end{pro}
\textsc{Proof}: Since $ B \geq 0 $ in $ \O$
and  using Lemma \ref{lem:3}
$$
\int_\O B v^{2k}\geq
\min_{\bar\O_{\a_k}}B \int_{\O_{\a_k}} v^{2k}
\geq \pi \min_{\bar\O_{\a_k}}B \intt \hv^{2k}
\geq \frac{1}{4^k} \min_{\bar\O_{\a_k}}B
\int_\O v^{2k}\ .
$$\qed

\begin{rem}
We stress that
estimate \equ{sti:fund} is not trivial
$($in the case $c_k(B)>0$$)$
since $ B $ could vanish   on $ \partial \O $
$($in particular, in this case,
 \equ{sti:fund} does not hold true in the whole
$ L^{2k} (\O)$$)$.
\end{rem}

\begin{rem}\label{ipdeb}
In light of Proposition \ref{pro:fund} we can prove Theorem 2
requiring only $H\geq 0$ in $\O$ and $H>0$ in $\overline{\O}_{\a_k}$,
instead of \equ{G:positiva}.
\end{rem}

In the following $\kappa_i$ will denote positive constants
depending only on $H$, $ \b $, $\mathcal{R}$, $k$
but not on $R$, $\e$. We also recall the notation
$ o(1)$ for a function tending to $0$ as $ \e \to 0 $
uniformly.

\subsection{The $L^{2k}$-estimate.}

Take $ \f:=\bv $ in the variational inequality
\equ{minore0}; $\f$ is an admissible variation since
$\langle \bv,\f  \rangle_{H^1}=\| \bv\|_{H^1}^2 >0 $.

\noindent
By
\equ{minore0},
$ \|w(\bv, \e )\|_E = O( \e ) $ (recall \equ{wC0*})
and \equ{R*u},
there exists  $0<\e_1\leq\e_0$ such
that
\begin{equation}\label{stima:L2k1}
\int_\O \b (\bv+H)^{2k}  \bv\leq 1\, \qquad\mbox{for}\quad
|\e|\leq \e_1\,,
\end{equation}
since
$$
\int_\O \b (\bv+H)^{2k}  \bv \leq \int_\O |{\cal R}_* \bv | +
\Big| \b (H + \bv + \bw )^{2k} - \b (H + \bv)^{2k} \Big| |\bv | \leq
o(1)\| \bv\|_{L^1}\leq o(1)R.
$$

\noindent
Noting that,  $\int_\O \b \bv^{2k+1}=0 $ by \equ{pari} with
$a(x,u) = \b(x) u^{2k}$ and $\b(\pi - x) = \b (x) $,
we derive
\begin{eqnarray}
  \int_\O\b (\bv+H)^{2k}  \bv &=& \int_\O\b
\Big[ (\bv+H)^{2k}-\bv^{2k} \Big] \bv \nonumber  \\
   &=& \int_\O 2k \b H \bv^{2k}+ \b
\sum_{j=0}^{2k-2}\left( {2k}\atop{j} \right)\bv^{j+1}H^{2k-j}
    \nonumber\\
   &\geq& 2k c_k(\b H) \| \bv\|_{L^{2k}}^{2k}-\kappa_1\| \bv\|_{L^{2k}}^{2k-1}-
   \kappa_2\| \bv\|_{L^{2k}}\ ,\label{sti:L2k1}
\end{eqnarray}
where $c_k(\b H)>0$ was defined in \equ{sti:costfund}
(recall \equ{bH})
and
 we have used Proposition \ref{pro:fund} and H\"{o}lder inequality
to estimate $\|\bv\|_{L^i} \leq C_{i,k} \|\bv\|_{L^{2k}}$ ($i\leq 2k-1$).

Finally, by \equ{stima:L2k1} and \equ{sti:L2k1} we deduce
\begin{equation}\label{sti:L2k_fin}
\| \bv\|_{L^{2k}}\leq \kappa_3
\, \qquad\mbox{for}\quad
|\e|\leq \e_1\,.
\end{equation}

\subsection{The $L^\infty$-estimate.}\label{subsec:Li}

To obtain the $ L^\infty $-estimate for $\bv $ we consider
an admissible variation $\f $,
which is a nonlinear function of $ \bv $,
and it is constructed as in [R67].

Let, for $M>0,$
\begin{equation}\label{q}
    q(\l):=q_M(\l):=\left\{
\begin{array}{ll}
0\ , & \mbox{if}\quad |\l|\leq M\\
\l-M\ , & \mbox{if}\quad \l\geq M\\
\l+M\ , & \mbox{if}\quad \l\leq M\ .
\end{array}
    \right.
\end{equation}
For $ \bv(t,x) = \bv_+(t,x)-\bv_-(t,x)=\hbv(t+x)-\hbv(t-x),$ we define
$$
\f := q_+-q_-:=q(\bv_+)-q(\bv_-) \in V \ .
$$
We take
$$
M := \frac12  \|\hbv\|_{L^\infty(\mathbb{T} )}
$$
and we can assume $ M > 0 $, i.e. $ \bv $ is not identically zero.

In [R67] it is proved that such $ \f $ is an admissible variation.
We report the proof for completeness.
By \equ{scalar}, it is sufficient
to prove that $\langle \bv, \f \rangle_{H^1} > 0$.

Using \equ{fubini2}, \equ{fubini6} and \equ{fubini7}
\begin{equation}
\langle \bv_+ - \bv_-, q ( \bv_+) - q ( \bv_-) \rangle_{H^1} =
\int_\O \bv_+ q ( \bv_+) + \bv_-q ( \bv_-)
 + 2 \int_\O q'(\bv_+  ) \Big[ (\bv_+)_t^2 + (\bv_+)_x^2  \Big]
\label{Rav}
\end{equation}
Since $ q $ is a monotone odd function of its argument and
by our choice of $ M $, $ \bv_\pm q ( \bv_\pm ) > 0 $ in
a positive measure set,
and,
since $ q' \geq 0 $, the second  term in
\equ{Rav} is non-negative.

We also have, since $ q $ is a monotone function,
\begin{equation}\label{vfpos}
\bv\f=(\bv_+-\bv_-)(q_+-q_-)=(\bv_+ - \bv_-)(q(\bv_+) - q(\bv_-))\geq 0\ .
\end{equation}

\noindent
Insert such $\f $ in the variational inequality \equ{minore0}.
Here the dominant term is $\int \b(\bv+H)^{2k}\f$,
in the sense that, by $\|w(v,\e) \|_E = O( \e )$ (recall \equ{wC0*}),
\equ{R*u} and $\|\bv\|_{L^\infty}\leq R$, we obtain
that there exists $0<\e_2\leq \e_1$ such that
\begin{equation}\label{stima:Li1}
\int_\O \b (\bv+H)^{2k}  \f\leq \| \f\|_{L^1}
\, \qquad\mbox{for}\quad
|\e|\leq \e_2\,.
\end{equation}

\noindent
Since $\int_\O \b \bv^{2k}\f=0$ by \equ{pari} and $\b (\pi - x)=
\b (x) $, we have (recall $ \b > 0 $)
\begin{equation}\label{sti:Lidom}
\int_\O \b(\bv+H)^{2k}\f =\int_\O \b \Big[ (\bv+H)^{2k}-\bv^{2k}\Big]\f
\geq \int_\O 2 k \b H \bv^{2k-1}\f-\kappa_4
\big( \|\bv\|_{L^\infty}^{2k-2}+1\big)
\| \f\|_{L^1}\,.
\end{equation}
We now estimate
the dominant term
$2 k \int_\O \b H \bv^{2k-1}\f$.
Since $\bv \f\geq 0$ and $\min_{\O_{1/4}} \b H>0$ (by \equ{bH})
\begin{equation}
   \int_\O 2k\bv^{2k-1}\b H\f = 2k \int_\O \b H (\bv\f)\bv^{2k-2}
   \geq \kappa_5 \int_{\O_{1/4}}  (\bv\f)\bv^{2k-2}\label{sti:Li1}\,.
\end{equation}
By \equ{stima:Li1}, \equ{sti:Lidom} and \equ{sti:Li1}, we have
\begin{equation}\label{sti:Li14}
\int_{\O_{1/4}}  \bv^{2k-1}\f\leq
\kappa_6 \big( \|\bv\|_{L^\infty(\O)}^{2k-2}+1\big)
\| \f\|_{L^1(\O)}\,.
\end{equation}
We have to give a lower bound of the positive integral
$\int_{\O_{1/4}}  \bv^{2k-1}\f=$
$\int_{\O_{1/4}}  (\bv\f)\bv^{2k-2}$
$=\int_{\O_{1/4}}  (\bv\f)(\bv_+ - \bv_-)^{2k-2}$.

\noindent
We first consider the  (more difficult) case $k\geq 2$,
in which the $L^{2k}$-estimate for $\bv $ obtained in the previous
subsection is needed, the (simpler)
case $k=1$ will be treated later.

Using
\equ{zialalletta}
we obtain
\begin{eqnarray}
  \int_{\O_{1/4}}  \bv^{2k-1}\f &\geq&
  \int_{\O_{1/4}}  \bv\f
\Big[ \bv_+^{2k-2}+\bv_-^{2k-2}-(2k-2)(\bv_+^{2k-3}\bv_-
+ \bv_+\bv_-^{2k-3}) \Big] \nonumber\\
   &=& 2\int_{\O_{1/4}}  \bv_+^{2k-1} q_+    - \bv_+^{2k-1} q_- + \bv_+^{2k-2} \bv_- q_- - \bv_+^{2k-2}\bv_- q_+
  \nonumber\\
   & & \qquad \ \ \ +(2k-2)[-
\bv_+^{2k-2}\bv_- q_+    + \bv_+^{2k-2}\bv_- q_- - \bv_+^{2k-3} \bv_-^2 q_- + \bv_+^{2k-3}\bv_-^2 q_+
        ]  \nonumber\\
   &\geq&
   2\int_{\O_{1/4}}  \bv_+^{2k-1} q_+    \label{sti:Li2} \\
   & &\!\!\!\!-2\int_{\O_{1/4}}      \bv_+^{2k-1} q_-  + (2k-1)\bv_+^{2k-2}\bv_- q_+
   +(2k-2) \bv_+^{2k-3} \bv_-^2 q_-
       \label{sti:Li3}
\end{eqnarray}
where in the  equality we have used \equ{fubini7} and in the last inequality
the fact that $\bv_+ q_+,$ $\bv_- q_-$ $\geq 0$ (since $\l q(\l)\geq 0$).

The dominant term is \equ{sti:Li2}.
Since $\l^{2k-1} q(\l)\geq M^{2k-1} |q(\l)|$, by
 \equ{fubini6} we obtain
\begin{equation}\label{v+q+}
2\int_{\O_{1/4}} \bv_+^{2k-1}q_+=2\pi(1-\frac24)\intt \hbv^{2k-1}(s)q(\hbv(s))ds\geq
\pi M^{2k-1} \| q(\hbv)\|_{L^1(\mathbb{T})}\ .
\end{equation}

\noindent
We now give an upper estimate of the three terms in \equ{sti:Li3}.
By \equ{fubini2}
\begin{eqnarray}
  \left| 2\int_{\O_{1/4}} \bv_+^{2k-1}q_- \right|
  &\leq&
  2\int_\O \left|  \bv_+^{2k-1} \right| |q_-| \leq
\| \hbv\|^{2k-1}_{L^{2k-1}(\mathbb{T})}
\| q(\hbv)\|_{L^1(\mathbb{T})}  \,, \nonumber
   \\
 \left| 2\int_{\O_{1/4}} \Big(\bv_+^{2k-2}q_+\Big)\bv_- \right|
  &\leq&
\| \hbv^{2k-2}q(\hbv)\|_{L^{1}(\mathbb{T})}
\|\hbv\|_{L^1(\mathbb{T})}
\leq (2M)^{2k-2}
\| q(\hbv)\|_{L^1(\mathbb{T})} \| \hbv\|_{L^1(\mathbb{T})}\nonumber
   \\
 \left| 2\int_{\O_{1/4}} \bv_+^{2k-3}\Big(\bv_-^2q_-\Big) \right|
   &\leq&
   \| \hbv^{2k-3}\|_{L^{1}(\mathbb{T})}
\|\hbv^2q(\hbv)\|_{L^1(\mathbb{T})}
\leq (2M)^2
\| \hbv\|^{2k-3}_{L^{2k-3}(\mathbb{T})}
\| q(\hbv)\|_{L^1(\mathbb{T})}\,.\nonumber
\end{eqnarray}
By the previous inequalities, \equ{v+q+}, H\"{o}lder
inequality\footnote{
To estimate
$\| \hbv\|_{L^j(\mathbb{T})} \leq C_{j,k}
\| \hbv\|_{L^{2k}(\mathbb{T})}$ for $j<2k.$}
and
\equ{sti:L2k_fin}, we finally have
\begin{equation}\label{sti:Liinterm}
    \int_{\O_{1/4}}  \bv^{2k-1}\f\geq
\pi M^{2k-1} \| q(\hbv)\|_{L^1(\mathbb{T})}-
\kappa_7 (M^{2k-2}+1) \| q(\hbv)\|_{L^1(\mathbb{T})}\,.
\end{equation}
Now we note that by \equ{fubini6}
\begin{equation}\label{fq}
    \int_\O |\f|\leq \int_\O |q(v_+)|+|q(v_-)|\
    =2\pi \| q(\hbv)\|_{L^1(\mathbb{T})}\,.
\end{equation}
We collect
\equ{sti:Li14}
and \equ{sti:Liinterm}
using \equ{fq} in order to obtain
\begin{equation}\label{pace}
M^{2k-1} \| q(\hbv)\|_{L^1(\mathbb{T})}\leq
\kappa_8
\Big(\|\bv \|_{L^\infty(\O)}^{2k-2} + M^{2k-2}+1\Big) \| q(\hbv)\|_{L^1(\mathbb{T})}\,.
\end{equation}
Since $ M:= \|\hbv\|_{L^\infty(\mathbb{T} )} \slash 2 $
hence $\|q(\hbv)\|_{L^\infty(\mathbb{T} )}= M$,
$\|\bv\|_{L^\infty(\O )}\leq$ $ 2\|\hbv\|_{L^\infty(\mathbb{T})}$
$= 4M$
and $\|\hbv\|_{L^\infty(\mathbb{T})}\neq 0$. Hence,
by \equ{pace},
$$
M^{2k-1} \| q(\hbv)\|_{L^1(\mathbb{T})}\leq
\kappa_9 \Big(M^{2k-2}+1\Big) \| q(\hbv)\|_{L^1(\mathbb{T})}
$$
and, dividing by $ \| q(\hbv)\|_{L^1(\mathbb{T})} \neq 0 $,
we finally obtain $M^{2k-1} \leq \kappa_9 (M^{2k-2}+1)$.

By our choice of $M$ the $L^\infty$-estimate follows
for $k\geq 2,$
\begin{equation}\label{sti:Li_fin}
    \| \bv\|_{L^\infty}\leq \kappa_{10}\, \qquad
\mbox{for}\quad |\e|\leq\e_2\,.
\end{equation}

We now briefly discuss the case $k=1$, which is simpler
and where a previous $L^2$-estimate for $\bv $
is not necessary to obtain \equ{sti:Li_fin}. In fact by \equ{stima:Li1}
and \equ{sti:Lidom} (with $k=1$), we obtain
\begin{equation}\label{sti:Licaso2}
\int_\O \b H\bv\f\leq \kappa_{11} \| \f\|_{L^1} \, .
\end{equation}
For $0 < \a<1/2$ to be chosen later, we have
\begin{equation}\label{sti:Lialpha}
\int_\O \b H\bv\f\geq \min_{\overline{\O}_\a} (\b H) \intOa \bv\f\,.
\end{equation}
We have to give a lower bound of
\begin{equation}\label{intOaq}
\intOa \bv\f=\intOa \bv_+q_+ +\bv_-q_- -\intOa \bv_+q_- +\bv_-q_+\,.
\end{equation}
By \equ{fubini6} and $\l q(\l)\geq M |q(\l)|$
\begin{equation}\label{pipi1}
\intOa \bv_+q_+=\intOa \bv_-q_-=\pi(1-2\a)\intt \hbv(s) q(\hbv(s))ds
\geq \pi(1-2\a) M \| q(\hbv)\|_{L^1(\mathbb{T})}\,.
\end{equation}
Moreover, since $\hbv$ has zero average, by \equ{fubini2}, we have
\begin{equation}\label{v-q+}
\left|\intOa \bv_+q_-\right|=\left|\intOa \bv_-q_+\right|
    \leq \frac12 \int_{-2\a\pi}^{2\a\pi}
   \intt |q(\hbv(y))\hbv(z+y)|dydz\leq 2\a\pi M
\| q(\hbv)\|_{L^1(\mathbb{T})}\,.
\end{equation}
Collecting \equ{intOaq}, \equ{pipi1} and \equ{v-q+} we obtain
\begin{equation}\label{pipi2}
\intOa \bv\f\geq 2\pi(1-6\a)M \| q(\hbv)\|_{L^1(\mathbb{T})}
\geq\frac12 \pi M \| q(\hbv)\|_{L^1(\mathbb{T})}\,,
\end{equation}
choosing $\a:=1/8.$
Collecting
\equ{sti:Licaso2},
\equ{sti:Lialpha}
and \equ{pipi2},
we obtain
$$
M \| q(\hbv)\|_{L^1(\mathbb{T})}\leq \kappa_{12} \|\f\|_{L^1(\O)}\,.
$$
Using \equ{fq} and dividing by $\| q(\hbv)\|_{L^1(\mathbb{T})} \neq 0 $
in the previous inequality we finally obtain
\equ{sti:Li_fin} also in the case $k=1.$

\medskip

\noindent

\subsection{The $H^1$-estimate.}\label{subsec:LH1}

We note that $\f:=-D_{-h}D_h \bv$ is an admissible variation, since
using \equ{perparti},
$$
\langle -D_{-h}D_h \bv, \bv \rangle_{H^1} =  \langle
D_h \bv, D_h \bv \rangle_{H^1} > 0.
$$

\noindent
Since $ \partial_t [\b (H+ \bv+\bw )^{2k}]$,
$ \partial_t \mathcal{R}_* \in L^2( \O) $ (see \equ{R*t}),
we have, as $ h \to 0$,
\begin{eqnarray*}
    &&\int_\O \b \big(H+ \bv+\bw\big)^{2k}\f
    \stackrel{\equ{perparti}}=\\
    &&\qquad=\int_\O \b
D_h\Big[\big(H+ \bv+\bw\big)^{2k}\Big]D_h \bv
    \stackrel{\equ{conL2quoz}}\longrightarrow
    \int_\O \b \Big[\big(H+ \bv+\bw\big)^{2k}\Big]_t \bv_t \ ,
\end{eqnarray*}
\begin{equation*}
\int_\O \mathcal{R}_*\f
\stackrel{\equ{perparti}}=
\int_\O D_h \mathcal{R}_* D_h \bv
\stackrel{\equ{conL2quoz}}\longrightarrow
\int_\O \partial_t \mathcal{R}_* \bv_t\,,
\end{equation*}
and, by the variational inequality \equ{minore0}, we obtain
\begin{equation}\label{minore0H1}
\int_\O \b \Big[(H+ \bv+\bw)^{2k}\Big]_t \bv_t
\leq
\int_\O - \partial_t \mathcal{R}_* \bv_t\,.
\end{equation}
By the $L^\infty$-estimate on  $\bv$
given in \equ{sti:Li_fin},
the Cauchy-Schwartz inequality
and \equ{R*t} we obtain
\begin{equation}\label{bla}
\left|
\int_\O  \partial_t \mathcal{R}_* \bv_t
\right|
\leq
\kappa_{13}
\| \bv_t \|_{L^2}+o(1)\| \bv_t \|_{L^2}^2 \ .
\end{equation}

\noindent
Since $\| \bw\|_E$ $=$ $\| w (\bv, \e)\|_E $ $= O(\e)$ (recall \equ{wC0*}),
again by  \equ{sti:Li_fin}and
the Cauchy-Schwartz inequality, we
find
\begin{eqnarray}
\int_\O \b \Big[\big(H+ \bv+\bw\big)^{2k}\Big]_t \bv_t
   &=&
2k\int_\O \b
\big(H+ \bv+\bw\big)^{2k-1} \big(H_t+ \bv_t+\bw_t\big)\bv_t\nonumber
   \\
   &\geq&
   2k \int_\O \b (\bv + H)^{2k-1}\bv_t^2
-o(1)\| \bv_t \|_{L^2}^2
   -\kappa_{14}
\| \bv_t \|_{L^2}\, .    \label{meriadoc}
\end{eqnarray}
Collecting \equ{minore0H1},\equ{bla} and \equ{meriadoc}, we obtain
\begin{equation}\label{pipino}
\int_\O \b (\bv + H)^{2k-1}\bv_t^2 \leq
\kappa_{15}
\| \bv_t \|_{L^2}+o(1)\| \bv_t \|_{L^2}^2\,.
\end{equation}
Since $\bv,\bv_t\in N$ and $\bv^{2k-1}\bv_t^2\in L^1(\O),$
it results $\int_\O \b \bv^{2k-1}\bv_t^2=0$ by \equ{pari2}. Using the
inequality \equ{zialalletta3} we obtain
\begin{eqnarray}
\int_\O \b (\bv + H)^{2k-1}\bv_t^2 &= &
\int_\O \b \Big[ (\bv + H)^{2k-1}-\bv^{2k-1} \Big] \bv_t^2 \nonumber \\
& \geq &
4^{1-k} \int_\O \b H^{2k-1}\bv_t^2
\geq  4^{1-k} c_1(\b H^{2k-1})\int_\O \bv_t^2 \label{sti:H1}
\end{eqnarray}
where $c_1(\cdot)$ was defined in \equ{sti:costfund} and $\b H^{2k-1}>0$ by \equ{bH}.
By \equ{pipino} and \equ{sti:H1} we get
$$
\| \bv_t \|_{L^2}^2
\leq
\kappa_{16}
\| \bv_t \|_{L^2}+o(1)\| \bv_t \|_{L^2}^2
$$
and we finally deduce that there exists a $0<\e_3\leq \e_2$ such that
\begin{equation}\label{sti:H1_fin1}
    \|\bv\|_{H^1}<\kappa_{17} \, \qquad \forall\, |\e|\leq \e_3\,.
\end{equation}
\medskip
\noindent
\textsc{{ Proof \ of \ Theorem \ 2-($i$) completed}}.
Defining $R:=\kappa_{17}$ and $\e_*:=\e_3$
we obtain, by \equ{sti:H1_fin1}, that
$$
\|\bv(\e)\|_{H^1} < R\qquad\qquad \forall\, |\e|\leq \e_*
$$
and $ \bv(\e) $ is an interior minimum or maximum point of $\Phi $
in $ B_{R} := \{ \|v\|_{H^1} < R  \} $. By Lemma \ref{natcon}
$\bu=$ $\bv+\bw=$ $ \bv(\e) + w (\bv(\e) , \e )$ is a weak solution
of \equ{equbis} and
\begin{equation}\label{u}
u:=\e \Big( H +  \bv(\e^{2k}) + w (\bv(\e^{2k}) , \e^{2k} ) \Big) \in E
\end{equation}
is a weak
solution of \equ{equ}-\equ{bc}-\equ{pc} satisfying
$\|u \|_E \leq C |\e |$. \qed

\begin{rem}\label{rem:non}
Let $u := u_\e = v_\e + w_\e $, $v_\e \in V$, $w_\e \in W $,
be a weak solution of $ \square u_\e = $ $ \e (g(x,u_\e) + h(t,x) ) $
where $g \in C( [0,\pi] \times \mathbb{R} ) $,
$ g(x,u) = g(x,-u) $ $=g(\pi - x,u) $. Suppose $u_\e $
satisfies $ \| u_\e \|_{L^\infty} \leq R $, $ \forall \e $ small.
We claim this implies $ h \in N^\bot $.
Indeed, $ w_\e $ satisfies the range equation
$ w_\e = $ $ \e \square^{-1} \Pi_{N^\bot} (g(x,u_\e) + h(t,x) ) $
and therefore
$\| w_\e \|_{L^\infty} \leq C |\e | $.
Moreover,
by the kernel equation
$ \Pi_N ( g(x, v_\e + w_\e ) + h(t,x)) = 0 $, and noting that
$ \Pi_N g(x, v_\e ) = 0 $ by \equ{pari}, we derive
$$
\Big|\Big| \Pi_N h(t,x) \Big|\Big|_{L^2} =
\Big|\Big|\Pi_N (g(x, v_\e + w_\e ) -g(x, v_\e )) \Big|\Big|_{L^2} \leq
\Big|\Big| g(x, v_\e + w_\e ) -g(x, v_\e ) \Big|\Big|_{L^2}
\to 0
$$
as $\| w_\e \|_{L^\infty} \to 0 $
because
$ g $ is uniformly continuous on $ [0,\pi] \times \{ |u| \leq C \}$.
Therefore $\Pi_N h = 0 $ and $h \in N^\bot $.
\end{rem}

\begin{rem}\label{multi} {{\bf (Multiplicity)}}
By \equ{pari},
any forcing term $ h(t,x) := - g (x, v_0(t,x) ) $,
$ v_0 \in V\setminus\{ 0\} $,
is in $ N^\bot $, if
$ g(x,u) = g(x,-u) $ $=g(\pi - x,u) $.
  Therefore the equation
$\square u = \e ( g(x,u) + h(t,x) )$
possesses, beyond the $\e$-small solution $u$ of Theorem 2,
also the other two (not small) solutions
$ \pm v_0 $.
\end{rem}

\subsection{Higher regularity and Classical Solutions}\label{sec:hig}

We now prove Theorem 2-(ii) obtaining more regularity
for the weak solution $u\in E$ of
\equ{equ}-\equ{bc}-\equ{pc} defined in \equ{u}.

Since $ \bv \in V := N \cap H^1 $ is
a critical point of $ \Phi : V \to \mathbb{R} $
\begin{equation}\label{psi=0}
    \int_\O
\big( \b (H+ \bu)^{2k} +
\mathcal{R}_* \big) \psi
=0\, \qquad \forall\, \psi \in N \cap H^1\,,
\end{equation}
which actually holds for any $\psi\in N$ since
$ N \cap H^1$ is  dense in $N$ with the $L^2$-topology\footnote{
Recall that $\big[\b (H+ \bu)^{2k} +
\mathcal{R}_*\big]\in L^2 (\O) $ since
$\b, H, \bu,\mathcal{R}_*$ are continuous functions.}.

 \noindent
Hence, taking $\psi:=\f_t$ for any $\f\in N \cap H^1$
in \equ{psi=0} and  integrating by parts, we find
\begin{eqnarray}
0 &=&
\int_\O
\partial_t  \Big( \b (H+ \bv+\bw)^{2k} +
\mathcal{R}_* \Big) \f\nonumber\\
&=&
\int_\O
\Big[
2k\b (H+ \bu)^{2k-1}(H_t+\bv_t+\bw_t)
+\partial_t\mathcal{R}_*
\Big]\f\qquad \label{critderbis}
\end{eqnarray}
for any $\f\in N \cap H^1$.
Since the term
into square brackets $[\ldots]$ in \equ{critderbis}
is in $L^2 (\O) $ then, again by the
$L^2$-density of $ N \cap H^1$   in $ N$,
\equ{critderbis}
actually holds for any $\f\in N$.

\noindent
Setting for brevity
$$
z:=z(t,x;\e):=\Big( t,x,\e^{\frac{1}{2k}}\big(H(t,x)+ \bu(t,x;\e)\big) \Big)\,,
$$
we can write, from \equ{R*}-\equ{R**},
$$
\mathcal{R}_*(t,x;\e)= \e^{-1}\mathcal{R}(z) \qquad {\rm and} \qquad
\partial_t\mathcal{R}_*=\e^{-1+\frac{1}{2k}}\partial_u \mathcal{R}(z)\bv_t+A
\,
$$
where
\begin{equation}\label{A}
    A(t,x;\e):=\e^{-1}\partial_t \mathcal{R}(z)
    +\e^{-1+\frac{1}{2k}}\partial_u \mathcal{R}(z)
    \big(H_t(t,x)+\bw_t(t,x;\e)\big) \ .
\end{equation}
Then
\equ{critderbis} becomes
\begin{equation}\label{critder}
    \int_\O
\Big[
2k\b (H+ \bu)^{2k-1}(H_t+\bv_t+\bw_t)
+
\e^{-1+\frac{1}{2k}}\partial_u \mathcal{R}(z)\bv_t+A
\Big]\f=0\qquad\ \ \ \forall\, \f\in N \ .
\end{equation}

For the remainder of this subsection we shall take $ \e \neq 0$,
and  $K_i$ will denote suitable positive constants possibly
depending\footnote{However, such $ K_i $ can be
taken \textsl{independently} of $\e$ if we assume the further
regularity hypothesis \equ{hbis} on ${\cal R}$,
see remark \ref{rem:Hp3} and remark \ref{rem:inde}.} also on $\e$.
\\[1mm]
\indent
Since we are assuming that $h\in H^j\cap C^{j-1}$, $j \geq 1 $,
then, by \equ{regdiG},
$H\in H^{j+1}\cap C^{j}$.
Hence, to prove that $u\in H^{j+1}\cap C^{j}$,
by \equ{u},
it is sufficient to show that $\bv,\bw\in$
$H^{j+1}\cap C^{j}$.

We first prove that

\begin{lem}\label{fire}
$ \bu \in C^1 (\oO) \cap H^2 ( \O ) $.
\end{lem}

\textsc{Proof}:
We shall divide the proof in three steps.
\\[1mm]
{\bf Step 1:}
$\bw \in C^1 ( \oO ) \cap H^2 ( \O )$ and
\begin{equation}\label{piccdibw}
\| \bw\|_{C^1}+\| \bw\|_{H^2}\leq K_1 \ | \e | \ .
\end{equation}

We have $\bu=$ $ \bv + \bw \in C( \ov{\O} ) \cap H^1 ( \O ) $,
 $ H \in H^{2}(\O) \cap
C^1( \ov{\O} ) $ and
$\b\in H^1\big((0,\pi)\big)$.
Moreover, since
$\mathcal{R}\in C^1(\oO\times \mathbb{R})$ and $z(\cdot;\e)\in C\cap H^1$,
then
$\mathcal{R}_*(\cdot;\e) \in C \cap H^1$ and
$\| \mathcal{R}_* (\cdot;\e)\|_{C} +
\| \mathcal{R}_*(\cdot;\e)\|_{H^1}\leq K_2$.
Hence
$f(t,x,\bu(t,x;\e);\e)$ $=\b(x)\big(H(t,x)+\bu(t,x;\e)\big)^{2k}$
$+\mathcal{R}_*(t,x;\e)$ $\in C\cap H^1$
and $\| f\|_{C}+\| f\|_{H^1}\leq K_3$.

Therefore, since
$ \bw  $ solves the range equation $\bw=\e \square^{-1}\Pi_{N^\bot} f$,
$\Pi_{N^\bot} $ satisfies \equ{Lovip}-\equ{proHj} and $ \square^{-1} $
satisfies \equ{regdiG}, we deduce that
$\bw \in C^1 ( \oO ) \cap $ $ H^2 ( \O )$ and
\equ{piccdibw}.
\\[2mm]
{\bf Step 2:} $ \bv_t \in L^\infty ( \O )$.
\\[1mm]
\indent
Let define
\begin{equation}\label{k18}
    \kappa_{18} :=
 k\, 4^{-k} \min_{\O_{1/8}}\big( \b H^{2k-1} \big)
\end{equation}
and $0<\e_4\leq \e_3$ such that, $\forall \, |\e| \leq \e_4 $
\begin{equation}\label{aragorn}
|\e|^{-1+\frac{1}{2k}}
\| \partial_u \mathcal{R}(z)\|_{L^\infty(\O)}\leq
\frac{\kappa_{18}}{2} \, , \
4 k \Big| \b (H+\bv + \bw )^{2k-1} - \b (H+\bv)^{2k-1} \Big|_{L^\infty}
\leq \pi \frac{\kappa_{18}}{2}
\end{equation}
(such $\e_4$ exists by \equ{h} and since $\|w\|_E = O(\e) $).
We claim that
\begin{eqnarray}\label{gimli}
    \int_\O \left[
    2k \b (H+\bu)^{2k-1}+ \e^{-1+\frac{1}{2k}} \partial_u \mathcal{R}(z)
    \right] v \f
    \geq
    \Big(10\pi\kappa_{18}M -\kappa_{19}\| \hv\|_{L^2(\mathbb{T})} \Big)
     \| q(\hv)\|_{L^1(\mathbb{T})}\ \
    \\
    \forall\, |\e|\leq \e_4\,, \quad
    \forall\, v=v_+ -v_- \in N\,,\ \  v_\pm(t,x)=\hv(t\pm x)\,,\quad
    \f:=q(v_+)-q(v_-)\in N\,,\nonumber
\end{eqnarray}
where $q=q_M$ ($ M > 0 $) was defined in \equ{q}.
Noting that $ \int_\O \b \bv^{2k-1} v \f = 0 $ by \equ{pari2},
$v \f \geq 0 $ and using \equ{zialalletta3},
\begin{eqnarray}
\int_\O 2k\b \big( H + \bv \big)^{2k-1} v \f &=&
\int_\O 2k\b \Big( \big( H + \bv \big)^{2k-1} - \bv^{2k-1}  \Big) v \f \geq
\int_\O 2k4^{1-k} \b H^{2k-1} v \f \nonumber\\
&\geq& 2k4^{1-k} \min_{\O_{1/8}}\big( \b H^{2k-1} \big)   \int_{\O_{1/8}} v \f
=8\kappa_{18}  \int_{\O_{1/8}} v \f\,.
\label{stiba}
\end{eqnarray}
Using \equ{fubini6}-\equ{fubini7} we obtain the
lower bound 
\begin{eqnarray}
\int_{\O_{1/8}} v \f & \geq & \int_{\O_{1/8}}
 v_+ q_+ + v_- q_- - \int_\O | v_- | |q_+ |
+ | v_+| | q_- | \nonumber \\
& \geq & 2 \pi \left(1- 2 \frac{1}{8} \right)
\int_0^{2\pi} \hv (s) q (\hv (s) )\, ds -
2 \int_{\O} |\hv (t-x) | | q ( \hv ( t + x )   )| \nonumber \\
& \geq &  \frac{3\pi}{2}
\int_0^{2\pi} \hv (s) q ( \hv (s) )\, ds -
\sqrt{2\pi}\| \hv \|_{L^2(\mathbb{T})} \| q ( \hv ) \|_{L^1(\mathbb{T})}\,,
\nonumber
\end{eqnarray}
and, by \equ{stiba}, we get
\begin{equation}\label{legolas}
\int_\O 2k\b \big( H + \bv \big)^{2k-1} v \f
\geq
12 \pi \kappa_{18}
\int_0^{2\pi} \hv (s) q ( \hv (s) )\, ds -
\kappa_{19}\| \hv \|_{L^2(\mathbb{T})} \| q ( \hv ) \|_{L^1(\mathbb{T})}    \,.
\end{equation}
Since $ \int_\O |v \f | = $ $ \int_\O v \f  = $
$ 2 \pi \intt \hv (s) q(\hv (s))\,ds$ and using
\equ{aragorn}, we get, $ \forall\, |\e|\leq \e_4 $,
 $$
\left |\int_\O
\e^{-1+\frac{1}{2k}} \partial_u \mathcal{R}(z)
     v \f
\right|
\leq
|\e|^{-1+\frac{1}{2k}}
\| \partial_u \mathcal{R}(z) \|_{L^\infty(\O)}\int_\O v\f
\leq \pi\kappa_{18} \intt \hv (s) q(\hv (s))\, ds \ .
$$
Therefore, using \equ{legolas} and \equ{aragorn}, we obtain
\begin{eqnarray}
&   & \int_\O \left[
    2k \b (H+\bu)^{2k-1}+ \e^{-1+\frac{1}{2k}} \partial_u \mathcal{R}(z)
    \right] v \f\nonumber\\
&  &  = \int_\O\!\!
    2k \b (H + \bv )^{2k-1} v \f
+ \e^{-1+\frac{1}{2k}} \partial_u \mathcal{R}(z) v \f
+ 2k \b \Big[ (H + \bv + \bw )^{2k-1} - (H + \bv )^{2k-1} \Big] v \f  \nonumber
 \\
&  & \geq
    10 \pi\kappa_{18}
\int_0^{2\pi} \hv (s) q ( \hv (s) )\, ds -
\kappa_{19}\| \hv \|_{L^2(\mathbb{T})} \| q ( \hv ) \|_{L^1(\mathbb{T})}    \,,
\qquad \forall\, |\e|\leq \e_4 \,. \label{balin}
\end{eqnarray}
Since $\l q(\l ) \geq M|q(\l )| $
$$
\intt \hv (s) q(\hv (s)) \,ds
\geq
M \intt |q(\hv (s))|\, ds  =
M\| q ( \hv ) \|_{L^1(\mathbb{T})}\,,
$$
and, by \equ{balin}, we finally get
\equ{gimli}.

\noindent
We now conclude the proof that $\bv_t\in L^\infty (\O ).$
Taking $ v:=\bv_t $ in \equ{gimli} and
$ \f:=q(\partial_t \bv_+)-$ $q(\partial_t \bv_-)$ we obtain\footnote{
Here $\bv(t,x)=\hbv(t+x)-\hbv(t-x)$ and so
$\bv_t(t,x)=\hbv'(t+x)-\hbv'(t-x)$.}
\begin{equation}\label{gloin}
    \int_\O \left[
    2k \b (H+\bu)^{2k-1}+ \e^{-1+\frac{1}{2k}} \partial_u \mathcal{R}(z)
    \right] \bv_t \f
    \geq
    \big(10\pi\kappa_{18}M -\kappa_{20} \big) \| q(\hbv')\|_{L^1(\mathbb{T})}\ \ \ \ \
    \forall\, |\e|\leq \e_4\,.
\end{equation}

\noindent
Note that, since $H,\bw\in C^1$,
by \equ{h} ($A$ is defined in \equ{A})
\begin{equation}\label{AA}
\| A\|_{L^\infty}\leq K_4\,.
\end{equation}
Finally
from \equ{critder}-\equ{piccdibw}-\equ{gloin}-\equ{AA}
we get
\begin{equation*}
\big(10\pi\kappa_{18}M -\kappa_{20} \big) \| q(\hbv')\|_{L^1(\mathbb{T})}
 \leq   K_5  \| \f \|_{L^1(\O)}\,,\qquad
\forall\, |\e|\leq \e_4 \,,
\end{equation*}
from which,
 recalling
$ \|\f \|_{L^1} \leq$
$ 2 \pi \| q ( \hat{\bv}' ) \|_{L^1(\mathbb{T})} $ (see \equ{fq}),
we have
\be\label{disf}
M\| q ( \hat{\bv}' ) \|_{L^1(\mathbb{T})}\leq
K_6 \| q ( \hat{\bv}' ) \|_{L^1(\mathbb{T})}\,.
\ee
We claim that \equ{disf}
implies $\hat{\bv}' \in L^\infty ( \mathbb{T} )$. Indeed, if
$\hat{\bv}' \not\in L^\infty ( \mathbb{T} )$, then for any $M > 0 $,
$\| q ( \hat{\bv}' ) \|_{L^1(\mathbb{T})} > 0 $ and \equ{disf} yields
$  M \leq K_6 $; hence  $\hat{\bv}' \in L^\infty ( \mathbb{T} )$.
Taking $M := \| \hat{\bv}' \|_{L^\infty(\mathbb{T})} / 2 $
we obtain, by \equ{disf}, $\|\hat{\bv}'\|_{L^\infty(\mathbb{T})} \leq  2 K_6$
and $\|\bv_t\|_{L^\infty(\O )} \leq 4 K_6 $ by \equ{isometriaL}.
\\[3mm]
{\bf Step 3: } $ \bv_t \in H^1 (\O)  $ (and hence $ \bv \in N \cap C^1 (\O)
\cap H^2 (\O )$).
\\[1mm]
\indent
We claim that
\begin{equation}\label{gimli2}
\int_\O\Big[
2k\b (H+ \bu)^{2k-1} +
\e^{-1+\frac{1}{2k}}\big(\partial_u \mathcal{R}(z)\big)\Big]
v^2
\geq
\kappa_{18}
\| v\|_{L^2}^2\,,\qquad\forall\, v\in N\,,\ \  |\e|\leq \e_4 \,,
\end{equation}
where $\kappa_{18}$ is defined in \equ{k18}
and $\e_4$ is defined in \equ{aragorn}.
Arguing as before, using
$ \int_\O \b \bv^{2k-1} v^2 = 0 $ (recall \equ{pari2}), \equ{zialalletta3}
and \equ{sti:fund}-\equ{sti:costfund}
\begin{eqnarray}
\int_\O 2k\b \big( H + \bv \big)^{2k-1} v^2 & = &
\int_\O 2k\b \Big( \big( H + \bv \big)^{2k-1} - \bv^{2k-1} \Big)
v^2 \nonumber \\
& \geq &
\int_\O 2k4^{1-k} \b H^{2k-1} v^2 \nonumber \\
& \geq &  2k4^{1-k} c_1( \b H^{2k-1}) \int_\O  v^2
= 2\kappa_{18} \| v \|_{L^2}^2\,.  \label{madai}
\end{eqnarray}
Since, by \equ{aragorn},
\begin{equation*}
  \left|
\int_\O
\e^{-1+\frac{1}{2k}}\big(\partial_u \mathcal{R}(z)\big)
v^2 \right| \leq
|\e|^{-1+\frac{1}{2k}}\|\partial_u \mathcal{R}(z)\|_{L^\infty}
\| v\|_{L^2}^2
\leq \frac{\kappa_{18}}{2} \| v\|_{L^2}^2\,,
\end{equation*}
and
\begin{eqnarray*}
\Big|
\int_\O 2 k \b [ (H + \bu )^{2k-1} - (H + \bv )^{2k-1} ] v^2
\Big| & \leq &
\Big| 2k\b [ (H + \bu )^{2k-1}- (H+\bv)^{2k-1} ]
\Big|_{L^\infty}  \| v\|_{L^2}^2  \\
& \leq &
\frac{\kappa_{18}}{2} \| v\|_{L^2}^2
\end{eqnarray*}
using \equ{madai} we prove \equ{gimli2}.

\noindent
Take $ \f = - D_{-h} D_h \bv_t \in N $ in \equ{critder}.
Integrating by part (recall \equ{perparti}
 and \equ{Leibniz})
equality \equ{critder},
we obtain
\begin{eqnarray}
0
&=&
\int_\O
D_h
\Big[
2k\b (H+ \bu)^{2k-1}(H_t+\bv_t+\bw_t)
+
\e^{-1+\frac{1}{2k}}\partial_u \mathcal{R}(z)\bv_t+A
\Big]
D_h\bv_t\nonumber\\
&=&
\int_\O\bigg[
2k\b (H+ \bu)^{2k-1} +
\e^{-1+\frac{1}{2k}}\big(\partial_u \mathcal{R}(z)\big)\bigg]
(D_h \bv_t)^2 \label{partedom}\\
& & \ \ +\bigg[
2k\b \Big(D_{h}\big( (H+\bu)^{2k-1} \big)\Big)\big(T_h\bv_t\big)
+
2k\b \Big(D_h\big( (H+\bu)^{2k-1} (H_t+\bw_t)\big)\Big)\nonumber
\\
& &\qquad\ \ \ +
\e^{-1+\frac{1}{2k}}\Big(D_{h}\big(\partial_u \mathcal{R}(z)\big)\Big)
\big(T_h \bv_t\big)
+D_h A
\bigg]D_h\bv_t\,.\nonumber
\end{eqnarray}
The dominant term here is \equ{partedom}.
Using \equ{gimli2} with $v:=D_h\bv_t,$ we get
\begin{equation}\label{partedom2}
    \int_\O\bigg[
2k\b (H+ \bu)^{2k-1} +
\e^{-1+\frac{1}{2k}}\big(\partial_u \mathcal{R}(z)\big)\bigg]
(D_h \bv_t)^2
\geq
\kappa_{18}
\| D_h \bv_t\|_{L^2}^2\,,\qquad
\forall\, |\e|\leq\e_4\,.
\end{equation}
We now estimate all the other terms.
Since $\| T_h\bv_t\|_{L^\infty (\O)}=$
$\| \bv_t\|_{L^\infty (\O) }\leq K_7$
and
\begin{equation}\label{bifur}
    \Big\| D_{h}\big( (H+\bu)^{2k-1} \big) \Big\|_{L^2}
\leq \Big\| \partial_t \big( (H+\bu)^{2k-1} \big) \Big\|_{L^2}
=(2k-1)
\Big\|  (H+\bu)^{2k-2} (H_t+\bu_t) \Big\|_{L^2}
\leq \kappa_{21}\,,
\end{equation}
we obtain
\begin{equation}\label{resti4}
 \left|
\int_\O\bigg[
2k\b \Big(D_{h}\big( (H+\bu)^{2k-1} \big)\Big)\big(T_h\bv_t\big)
\bigg]
D_h \bv_t
\right|
\leq
K_8
\| D_h \bv_t\|_{L^2}
\,.
\end{equation}
Since
$H\in C^1\cap H^2,$
$\bv_t\in L^\infty$, $\bw\in C^1\cap H^2,$
we have
\begin{eqnarray*}
& \Big\| D_h\big( (H+\bu)^{2k-1} (H_t+\bw_t)\big) \Big\|_{L^2}
   \leq
\Big\| \partial_t \big( (H+\bu)^{2k-1} (H_t+\bw_t)\big) \Big\|_{L^2}
   \qquad\qquad\qquad\qquad\\
   &=
\Big\| (2k-1) (H+\bu)^{2k-2} (H_t+\bu_t) (H_t+\bw_t)+
(H+\bu)^{2k-1} (H_{tt}+\bw_{tt}) \Big\|_{L^2}
   \leq
   K_9\,,
\end{eqnarray*}
and we deduce
\begin{equation}\label{resti5}
 \left|
\int_\O
\bigg[
2k\b
 \Big(D_h\big( (H+\bu)^{2k-1} (H_t+\bw_t)\big)\Big)
\bigg]
D_h \bv_t
\right|
\leq
K_9
\| D_h \bv_t\|_{L^2}
\,.
\end{equation}
From
$$
\e^{-1+\frac{1}{2k}}
\partial_t\big(\partial_u \mathcal{R}(z)\big)
=
\e^{-1+\frac{1}{2k}}
\partial_{tu}^2 \mathcal{R}(z)
+
\e^{-1+\frac{1}{k}}
\partial_{uu}^2 \mathcal{R}(z) (H_t+\bu_t)
$$
we derive
\begin{equation}\label{bombur}
\Big\| \e^{-1+\frac{1}{2k}}
D_h \big(\partial_u \mathcal{R}(z)\big) \Big\|_{L^2}
\leq
\Big\| \e^{-1+\frac{1}{2k}}
\partial_t \big(\partial_u \mathcal{R}(z)\big)\Big\|_{L^2}
\leq K_{10}
\end{equation}
and  (recall $\bv_t\in L^\infty$)
\begin{equation}\label{resti6}
 \left|
\int_\O
\bigg[
\e^{-1+\frac{1}{2k}}\Big(D_{h}\big(\partial_u \mathcal{R}(z)\big)\Big)\bv_t
\bigg]
D_h \bv_t
\right|
\leq
K_{10} \| \bv_t\| \
\| D_h \bv_t\|_{L^2} = K_{11} \| \bv_t\|
\,.
\end{equation}
We finally estimate the term
$$
\left|
\int_\O
(D_h A )(D_h\bv_t)
\right|
\leq
\| D_h A\|_{L^2} \| D_h \bv_t\|_{L^2}
\stackrel{\equ{sti:diffquo}}\leq
\| \partial_t  A\|_{L^2} \| D_h \bv_t\|_{L^2}\, .
$$
($ A $ is defined in \equ{A}). Since
\begin{eqnarray*}
 \partial_t A &=&
\e^{-1}
\partial_{tt}^2 \mathcal{R} (z)+
 \e^{-1+\frac{1}{2k}}
\partial_{tu}^2 \mathcal{R} (z)(2H_t+\bv_t+2\bw_t) \\
   & &+
\e^{-1+\frac{1}{k}}
\partial_{uu}^2 \mathcal{R} (z)(H_t+\bu_t)(H_t+\bw_t)+
 \e^{-1+\frac{1}{2k}}
\partial_{u} \mathcal{R} (z)(H_{tt}+\bw_{tt})\, ,
\end{eqnarray*}
by \equ{h}, 
and using that $\bv_t \in L^\infty (\O)$,
 $\bw,H \in H^2 (\O) \cap C^1(\O)$, then
 $\|  \partial_t A \|_{L^2} \leq K_{12}$ and
\begin{equation}\label{resti7}
\left|
\int_\O
(D_h A )(D_h\bv_t)
\right|
\leq K_{12}\| D_h \bv_t\|_{L^2}\,.
\end{equation}
Recollecting
\equ{partedom}, \equ{partedom2},
\equ{resti4}, \equ{resti5}, \equ{resti6} and \equ{resti7}
we obtain
$$
\kappa_{18} \|D_h \bv_t\|_{L^2}^2
\leq  K_{13}  \| D_h \bv_t\|_{L^2(\O)}
$$
and
$ \| D_h \bv_t \|_{L^2( \O )} \leq $ $ K_{14} $ for all
$ h $. Therefore, by \equ{weakd}, $ \bv_t \in N \cap H^1 $ and
$ \| \bv_{tt} \|_{L^2( \O )} \leq $ $ K_{14} $.\qed

\medskip

We now prove Theorem 2-(ii)
by induction over $j \geq 1 $.

\begin{lem}
Assuming
$ \bv$, $ \bw \in C^{j-1}( \ov{\O} ) \cap H^{j} ( \O ) $, then
$ \bv$, $ \bw \in C^{j}( \ov{\O} ) \cap H^{j+1} ( \O ) $.
\end{lem}

\textsc{Proof}: Again we divide the proof in three steps.
\\[1mm]
{\bf Step 1:} $ \bw \in $ $C^{j}( \ov{\O} ) \cap H^{j+1} ( \O ) $
and
\begin{equation}\label{piccdibwj}
\| \bw\|_{C^j}+\| \bw\|_{H^{j+1}}\leq K_1^{(j)} | \e | \,.
\end{equation}

By hypotheses
$ H \in $ $C^{j}( \ov{\O} ) \cap$ $ H^{j+1} ( \O ) $
(since $h\in$ $C^{j-1}( \ov{\O} ) \cap H^{j} ( \O ) $),
$\b\in H^j\big( (0,\pi)\big)$,
 $\mathcal{R}\in C^j(\oO\times \mathbb{R})$ and $z(\cdot;\e)\in $
 $ C^{j-1}( \ov{\O} ) \cap H^{j} ( \O ) $.
Hence
$\mathcal{R}_*(\cdot;\e) \in $ $C^{j-1}( \ov{\O} ) \cap H^{j} ( \O ) $ and
$\| \mathcal{R}_* (\cdot;\e)\|_{C^{j-1}}+\|
\mathcal{R}_*(\cdot;\e)\|_{H^{j}}\leq K_2^{(j)}$.
Hence
$ f (t,x,\bu(t,x;\e);\e)$
$ = \b(x)\big(H(t,x)+\bu(t,x;\e)\big)^{2k}$
$+\mathcal{R}_*(t,x;\e)$ $\in C^{j-1}( \ov{\O} ) \cap H^{j} ( \O ) $
and $ \| f\|_{C^{j-1}}+\| f\|_{H^{j}}\leq K_3^{(j)}$.

Since
$ \bw  $ solves the range equation $ \bw=\e \square^{-1}\Pi_{N^\bot} f$,
$\Pi_{N^\bot} $ satisfies \equ{Lovip}-\equ{proHj} and $ \square^{-1} $
satisfies \equ{regdiG}, we conclude that
$ \bw \in $ $C^{j}( \ov{\O} ) \cap H^{j+1} ( \O ) $ and that
\equ{piccdibwj} holds.
\\[2mm]
{\bf Step 2:}  $\partial_t^j\bv\in L^\infty (\O ) $.
\\[2mm]
\indent
Reasoning as for \equ{critderbis}, we get
\begin{eqnarray}
0 &=&
\int_\O
\partial_t^j  \Big( \b (H+ \bv+\bw)^{2k} +
\mathcal{R}_* \Big) \f\nonumber\\
&=&
\int_\O
\Big[
2k\b (H+ \bu)^{2k-1}
+
\e^{-1+\frac{1}{2k}}\partial_u \mathcal{R}(z)
\Big](\partial_t^j\bv)\f
+\mathcal{F}^{(j)}\f
 \label{critderbisj}
\end{eqnarray}
for any $\f\in N$.
Here
$\mathcal{F}^{(j)}$
depends polynomially on $k,\b$ and
$$
    \e^{-1+\frac{n}{2k}}\partial_t^l\partial_u^n\mathcal{R}(z)\
\ \ \ l+n= j\,,
\qquad
\partial_t^l H\,, \ \partial_t^l \bw\ \ \ \ l\leq j\,,
 \qquad \partial_t^l \bv\ \ \ \ l\leq j-1\,,
$$
but \textsl{not} on $\partial_t^j \bv.$

\noindent
Choosing in \equ{gimli} $ v = \partial_t^j \bv  $ and so
$\f:=q(\partial_t^j \bv_+)-q(\partial_t^j \bv_-)$,
we get (recall that $\bv\in H^j$),
$\forall\, |\e|\leq \e_4,$
\begin{equation}\label{gloinj}
\int_\O
\Big[
2k\b (H+ \bu)^{2k-1}
+
\e^{-1+\frac{1}{2k}}\partial_u \mathcal{R}(z)
\Big](\partial_t^j\bv)\f
\geq
\bigg( 10\pi\kappa_{18}M-\kappa_{24}\bigg)
\| q(\hbv^{(j)})\|_{L^1(\mathbb{T})}\,,
\end{equation}
where
$\hbv^{(j)}(\xi):=\frac{d^j}{d\xi^j}\hbv.$
On the other hand, by \equ{critderbisj}
we get
\begin{equation}\label{AAj}
\int_\O
\Big[
2k\b (H+ \bu)^{2k-1}
+
\e^{-1+\frac{1}{2k}}\partial_u \mathcal{R}(z)
\Big](\partial_t^j\bv)\f
\leq
\| \mathcal{F}^{(j)}\|_{L^\infty}\|\f\|_{L^1}
\leq
K_4^{(j)}\|\f\|_{L^1}\,.
\end{equation}
Since
$ \|\f \|_{L^1} \leq$
$ 2 \pi \| q ( \hbv^{(j)} ) \|_{L^1(\mathbb{T})} $ (see \equ{fq}),
from \equ{gloinj} and \equ{AAj} we get
\begin{equation*}
\big(10\pi\kappa_{18}M -\kappa_{24} \big)
\| q(\hbv^{(j)})\|_{L^1(\mathbb{T})}
 \leq   K_5^{(j)} \| q(\hbv^{(j)})\|_{L^1(\mathbb{T})}\,,\qquad
\forall\, |\e|\leq \e_4\,.
\end{equation*}
Then
$$
M\| q ( \hbv^{(j)} ) \|_{L^1(\mathbb{T})}\leq
K_6^{(j)} \| q ( \hbv^{(j)} ) \|_{L^1(\mathbb{T})}\,,
\qquad
\forall\, |\e|\leq \e_4\,.
$$
Arguing as in  \equ{disf}
we get $ \hbv^{(j)}\in L^\infty ( \mathbb{T} )$.
Finally
$\partial_t^j\bv\in L^\infty (\O ) $
and
$\|\partial_t^j\bv\|_{L^\infty (\O )}\leq 4 K_6^{(j)}.$

\medskip

\noindent
{\bf Step 3:} We now prove that
$\partial_t^j\bv\in H^1  $ (and hence $\bv\in N\cap C^j \cap H^{j+1}$).
\\[2mm]
\indent
Choosing $\f:=-D_{-h}D_h \partial_t^j\bv$ in
\equ{critderbisj}, integrating by parts
(recall \equ{perparti} and \equ{Leibniz})
\begin{eqnarray}
0
&=&
\int_\O
\Big[
2k\b (H+ \bu)^{2k-1}
+
\e^{-1+\frac{1}{2k}}\partial_u \mathcal{R}(z)
\Big]\big(D_h(\partial_t^j\bv)\big)^2
\label{critderbisj2}\\
& &
+
\bigg[
\Big(
D_h
\big(
2k\b (H+ \bu)^{2k-1}
+
\e^{-1+\frac{1}{2k}}\partial_u \mathcal{R}(z)
\big)\Big)
\big(T_h(\partial_t^j\bv)\big)
+
\big(D_h (\mathcal{F}^{(j)})\big)
\bigg]
\big(D_h(\partial_t^j\bv)\big)\,.
\nonumber
\end{eqnarray}
Using \equ{gimli2} we get
\begin{equation}\label{dori}
\int_\O
\Big[
2k\b (H+ \bu)^{2k-1}
+
\e^{-1+\frac{1}{2k}}\partial_u \mathcal{R}(z)
\Big]\big(D_h(\partial_t^j\bv)\big)^2
\geq
\kappa_{18}
\Big\| D_h(\partial_t^j\bv) \Big\|_{L^2}^2\,,\qquad
\forall\, |\e|\leq \e_4\,.
\end{equation}
From \equ{bifur}, \equ{bombur} and since $ \partial_t^j \bv
\in L^\infty (\O ) $,
\begin{eqnarray}
  \left|
\int_\O
\Big(
D_h
\big(
2k\b (H+ \bu)^{2k-1}
+
\e^{-1+\frac{1}{2k}}\partial_u \mathcal{R}(z)
\big)\Big)
\big(T_h(\partial_t^j\bv)\big)
\big(D_h(\partial_t^j\bv)\big)
  \right|
  \nonumber\\
\leq  4 K_6^{(j)}
\Big\| D_h \big( 2k\b (H+ \bu)^{2k-1} +
\e^{-1+\frac{1}{2k}}\partial_u \mathcal{R}(z)
\big) \Big\|_{L^2}
\| D_h(\partial_t^j\bv) \|_{L^2}
\nonumber\\
\leq K_7^{(j)}
\| D_h(\partial_t^j\bv) \|_{L^2}
\label{K7j}\,.
\end{eqnarray}
We note that
$\partial_t \mathcal{F}^{(j)}$
is a polynomial in $k, \b, $
$$
    \e^{-1+\frac{n}{2k}}\partial_t^l\partial_u^n\mathcal{R}(z)\
\ \ \ l+n= j+1\,,
\qquad
\partial_t^l H\,, \ \partial_t^l \bw\ \ \ \ l\leq j+1\,,
 \qquad \partial_t^l \bv\ \ \ \ l\leq j\,,
$$
and that
the terms $ \partial_t^{j+1} H $, $ \partial_t^{j+1} \bw \in L^2 (\O )$
(recall that $\bw, H \in H^{j+1} (\O)$ by Step $1$)
appear only linearly (with no powers). Hence,
using that $ \partial_t^l H $,  $\partial_t^l \bw $,
$\partial_t^l \bv \in L^{\infty}( \O ) $, $ \forall l \leq j $,
\begin{eqnarray}
\left|\int_\O   \big(D_h (\mathcal{F}^{(j)})\big)
\big(D_h(\partial_t^j\bv)\big)
\right|
&\leq&
\| D_h  (\mathcal{F}^{(j)})\|_{L^2}
\| D_h(\partial_t^j\bv) \|_{L^2}
\nonumber\\
&\stackrel{\equ{sti:diffquo}}
\leq&
\| \partial_t (\mathcal{F}^{(j)})\|_{L^2}
\| D_h(\partial_t^j\bv) \|_{L^2}
\nonumber\\
&\leq&
K_8^{(j)}
\| D_h(\partial_t^j\bv) \|_{L^2} \ .
\label{K8j}
\end{eqnarray}
Finally, by
\equ{critderbisj2}-\equ{dori}-\equ{K7j}-\equ{K8j}
we get
$$
\kappa_{18}
\| D_h(\partial_t^j\bv) \|_{L^2}^2
\leq
K_9^{(j)}
\| D_h(\partial_t^j\bv) \|_{L^2},
$$
and therefore
$  \| D_h(\partial_t^j\bv) \|_{L^2} \leq $$ K_{10}^{(j)} $.
By \equ{weakd}, we conclude the proof  obtaining
$\partial_t^{j+1} \bv \in $ $ L^2 (\O )$
and  $ \| \partial_t^{j+1}\bv \|_{L^2} \leq $
$ K_{10}^{(j)} $. \qed

\begin{rem}\label{rem:inde}
If $H$, $\bu \in H^i ( \O) \cap C^{i-1} (\oO ) $ $(0 \leq i \leq j) $
and \equ{hbis} holds, then $R_* $ and $\partial_t R_* $ are
bounded in $ H^i (\O ) \cap C^{i-1} (\oO ) $
by some constant $\kappa_i $ independent
of $\e $. In this case, the constants $ K_i $ of this section
can be taken  independently of $\e $, obtaining the estimates
\equ{sti:regol}.
\end{rem}

\subsection{Proof of Theorem 1}

The following Proposition
is a sort of ``maximum principle" for the
wave equation \equ{equ}-\equ{bc}-\equ{pc}.

\begin{pro}\label{lem:6}
Let $h\in N^\bot$, $h>0$ (or $h\geq 0$) a.e. in $\O.$
Then there exists
a weak solution $H\in E$ 
of $\square H=h$ satisfying $H>0$ (or $ H \geq 0 $).
In particular we can choose
\begin{equation}
H(t,x)   :=
\frac12 \int_0^{\kappa} \int_{t-x-\xi}^{t-x+\xi} h(\tau,\xi) \, d\tau\, d\xi
-\frac12 \int_{\kappa}^x \int_{t-x+\xi}^{t+x-\xi} h(\tau,\xi) \, d\tau\, d\xi
\label{Hc1}
\end{equation}
for a suitable ${\kappa}\in(0,\pi)$.
Moreover $h\in C^{j-1} \Rightarrow H\in C^{j}$ and
$h\in H^{j} \Rightarrow H\in H^{j+1}$, for $j\geq 1$.
\end{pro}
\textsc{Proof}:
We consider the case $h>0$, the case $h\geq 0$ being similar.
\\[1mm]
\noindent
{\bf Step 1:} {\sl $ H $ defined in \equ{Hc1} belongs to
$  H^1 (\O) \cap C^{1\slash 2}(\oO) $
for any ${\kappa}\in (0,\pi)$ and
\begin{eqnarray}
  2 (\partial_t H)(t,x) & = &
  \int_0^\kappa \Big( h(t-x+\xi,\xi)-h(t-x-\xi,\xi) \Big)d\xi \nonumber\\
  & & -\int_\kappa^x \Big( h(t+x-\xi,\xi) -
h(t-x+\xi,\xi) \Big)d\xi\label{erika} \in L^2 (\O ) \ ,  \\
  2 (\partial_x H)(t,x) & = &
  \int_0^\kappa \Big(- h(t-x+\xi,\xi) + h(t-x-\xi,\xi) \Big)d\xi \nonumber\\
  & & -\int_\kappa^x \Big( h(t+x-\xi,\xi) +
h(t-x+\xi,\xi) \Big)d\xi \in L^2 (\O ) \ .\label{ale}
\end{eqnarray} }

\noindent
We shall prove that the first addendum in the r.h.s. of
\equ{Hc1}
$$
H_1(t,x):=\frac12 \int_0^{\kappa}
\int_{t-x-\xi}^{t-x+\xi} h(\tau,\xi) \, d\tau\, d\xi
$$
belongs to $ C^{1\slash 2}(\oO) \cap H^1 (\O) $,
the second addendum being analogous. Defining
$$
T(t,x):=T(t,x;{\kappa}):=\left\{
(\tau,\xi)\in\O\ \ \ \Big|\ \ \ t-x-\xi< \tau< t-x+\xi\,,\ \ 0<\xi< {\kappa}
 \right\}
$$
we can write
$ H_1(t,x):= $ $ (1\slash 2) \int_{T(t,x)} h(\tau,\xi) \, d\tau\, d\xi $.

Since $\mbox{meas}\big( T(t,x;{\kappa}) \big)={\kappa}^2\leq \pi^2$
we derive that $H_1$ is uniformly bounded by
$$
\big| H_1(t,x) \big|\leq\frac12
\int_\O {\bf{1}}_{T(t,x)}(\tau,\xi)|h(\tau,\xi)| \, d\tau\, d\xi
\leq \frac{\pi}{2}\| h\|_{L^2(\O)}
$$
using Cauchy-Schwartz inequality.

For $i=1,2$ and   $(t_i,x_i)$ $\in \O,$ let define
$T_i:=$$T(t_i,x_i)$. It results
$$
\mbox{meas}(T_1\setminus T_2)=
\mbox{meas}(T_2\setminus T_1)
\leq
\pi \big(
|t_1-t_2|+|x_1-x_2|
\big)
$$
and, using again Cauchy-Schwartz inequality,
$$
\left|
H_1(t_1,x_1)-H_1(t_2,x_2)
\right|
\leq
\frac12 \int_{T_1\setminus T_2}|h| +
\frac12 \int_{T_2\setminus T_1}|h|
\leq
\sqrt{\pi}
\big(
|t_1-t_2|+|x_1-x_2|
\big)^{1/2}
\| h\|_{L^2(\O)} \ .
$$
Therefore we have proved $ H_1\in C^{1/2}(\oO) $.

\noindent
We now prove that $H_1\in H^1 (\O)$ and that
$\partial_t H_1=$ $-\partial_x H_1 = $$ f_1 $ where
$$
f_1(t,x) :=
\frac12 \int_0^\kappa \Big( h(t-x+\xi,\xi)-h(t-x-\xi,\xi)\Big) \, d\xi
\in L^2(\O) \ .
$$
We first justify that $ f_1 \in L^2(\O) $. Since
$$
|h(t-x+\xi,\xi)|^2+|h(t-x-\xi,\xi)|^2
\geq \frac12
|h(t-x+\xi,\xi)-h(t-x-\xi,\xi)|^2\,,
$$
by periodicity w.r.t. $t$ we obtain that, $\forall\, x\in(0,\pi)$,
$$
\| h\|^2_{L^2(\O)}
\geq
\int_0^{\kappa}\intt |h(t,\xi)|^2dt\, d\xi
\geq
\frac14
\int_0^{\kappa}\intt |h(t-x+\xi,\xi)-h(t-x-\xi,\xi)|^2dt\, d\xi\,.
$$
Integrating the previous inequality in the variable $x$ between $0$ and $\pi$,
applying Fubini Theorem and Cauchy-Schwartz inequality, we deduce
\begin{eqnarray*}
   \pi \| h\|^2_{L^2(\O)}
   &\geq&
   \frac14 \int_0^\pi  \int_0^{\kappa}\intt
|h(t-x+\xi,\xi)-h(t-x-\xi,\xi)|^2\, dt\, d\xi\, dx  \\
   &=&  \frac14 \int_\O \int_0^\kappa
|h(t-x+\xi,\xi)-h(t-x-\xi,\xi)|^2\,d\xi\, dt\, dx \ \\
   &\geq &  \frac{1}{4\kappa}
\int_\O \left( \int_0^\kappa
|h(t-x+\xi,\xi)-h(t-x-\xi,\xi)|\,d\xi\right)^2\, dt\, dx\, .  \\
&\geq &  \frac{1}{\kappa}
\int_\O f_1^2 (t,x)\, dt\, dx = \frac{1}{\kappa} \| f_1 \|_{L^2}^2 \ .
\end{eqnarray*}
Finally, we prove that $ \partial_x H_1 = - f_1$,
being $ \partial_t H_1 =$ $f_1$  analogous.
By Fubini Theorem
\begin{eqnarray*}
  2 \int_\O H_1 \f_x
  &=&
  \intt \int_0^\kappa \int_0^\pi \int_{t-x-\xi}^{t-x+\xi}
  h(\tau,\xi)\f_x (t,x)\, d\tau\, dx\, d\xi\, dt \\
   &=& \intt \int_0^\kappa \bigg( \int_{t+\xi-\pi}^{t+\xi} h
   \int_0^{t+\xi-\tau} \f_x \, dx\, d\tau
   + \int_{t-\xi}^{t+\xi-\pi} h
    \int_0^\pi \f_x \, dx\, d\tau \\
   & & \qquad \qquad\qquad\qquad\qquad\qquad\qquad\ \ \ \ +
    \int_{t-\xi-\pi}^{t-\xi} h
    \int_{t-\xi-\tau}^\pi \f_x \, dx\, d\tau\bigg)\, d\xi\, dt \\
   &=&
   \intt \int_0^\kappa \bigg( \int_{t+\xi-\pi}^{t+\xi} h(\tau,\xi)
   \f(t,t+\xi-\tau) \, d\tau\\
   & &\qquad \qquad\qquad-
   \int_{t-\xi-\pi}^{t-\xi} h(\tau,\xi)
   \f(t,t-\xi-\tau) \, d\tau
   \bigg)\, d\xi\, dt  \\
   &=&
   \intt \int_0^\kappa \bigg( -\int_\pi^0 h(t+\xi-x,\xi) \f(t,x)
    \, dx\\
   & &\qquad \qquad\qquad +
\int_\pi^0 h(t-\xi-x,\xi) \f(t,x)
    \, dx
   \bigg)\, d\xi\, dt  =  2 \int_\O f_1 \f\,.
\end{eqnarray*}
With analogue computations for the second addendum of $ H $
in \equ{Hc1} we derive \equ{erika} and \equ{ale}.
\\[1mm]
{\bf Step 2:}
{\sl There exists $ { \kappa } \in (0,\pi) $ such that $ H(t,x) $
verifies the Dirichlet boundary conditions
$H(t,0)=$ $H(t,\pi) = 0 $ $\forall\, t\in\mathbb{T} $.}
\\[1mm]
\indent
By \equ{Hc1}, the function $ H $ satisfies, for any $ \kappa \in (0, \pi) $,
$ H(t,0) = 0 $ $ \forall \, t \in \mathbb{T} $. It remains
to find $ \kappa $ imposing $ H (t,\pi) = 0 $.
Taking $x=\pi$ in \equ{Hc1} we obtain
\begin{eqnarray}
    H(t,\pi)
    &:=&
\frac12 \int_0^{\kappa} \int_{t-\pi-\xi}^{t-\pi+\xi} h(\tau,\xi) \, d\tau\, d\xi
-\frac12 \int_{\kappa}^\pi \int_{t-\pi+\xi}^{t+\pi-\xi} h(\tau,\xi) \, d\tau\, d\xi
\nonumber\\
   &=&
\frac12 \int_0^\pi \int_{t-\pi-\xi}^{t-\pi+\xi} h(\tau,\xi) \, d\tau\, d\xi
-\frac12 \int_{\kappa}^\pi \int_{t-\pi-\xi}^{t+\pi-\xi} h(\tau,\xi) \, d\tau\, d\xi
\nonumber\\
   &=&
   \frac12 \int_0^\pi \int_{-\xi}^{\xi} h(\tau,\xi) \, d\tau\, d\xi
   -\frac12 \int_{\kappa}^\pi \int_{-\pi}^{\pi} h(\tau,\xi) \, d\tau\, d\xi
     \ =: \ c -\chi({\kappa})\,,
   \label{carbonara}
\end{eqnarray}
where in the last line we have used the periodicity of
$h(\cdot,\xi)$ and \equ{clovi}.

In order to prove that $ H( t, \pi) = 0 $, $ \forall t \in \mathbb{T} $,
we need only to solve $ \chi({\kappa}) = c $.
By the absolute continuity of the integral
(with respect to the two-dimensional measure $d\tau\, d\xi$)
$ \chi({\kappa})$ is a 
continuous function.
Moreover, since $h>0$ a.e. in $\O$,
$\chi(0)> c >0$. Finally $\chi(\pi)=0$ and
therefore, by continuity,
there exists ${\kappa}\in(0,\pi)$ solving
 $\chi({\kappa})= c $.
\\[1mm]
{\bf Step 3:} {\sl $ H \in E $ is a weak solution of $ \square H = h $,
namely
\begin{equation}\label{raffi}
\int_\O \f_t H_t-\f_x H_x +\f h = 0 \,,
\qquad \forall\, \f\in C^1_0(\oO)\,.
\end{equation}}
By Fubini Theorem and periodicity  we get
\begin{eqnarray}
   & &  \int_\O \bigg(
   \f_t(t,x)\int_0^\kappa  h(t-x+\xi,\xi)d\xi
   + \f_x(t,x)\int_0^\kappa  h(t-x+\xi,\xi)d\xi
   \bigg)\,dt\, dx \nonumber\\
   & & =
   \int_0^\pi \int_0^\kappa \intt  \Big(
   \f_t(t,x)  h(t-x+\xi,\xi)
   + \f_x(t,x)  h(t-x+\xi,\xi)
   \Big)\,dt\,d\xi\, dx
   \nonumber\\
   & &=
   \int_0^\pi \int_0^\kappa \intt  \Big(
   \f_t(t+x-\xi,x)+
   \f_x(t+x-\xi,x)
   \Big)h(t,\xi)\,dt\,d\xi\, dx
   \nonumber\\
   & & =\int_0^\kappa\intt h(t,\xi) \int_0^\pi
   \frac{d}{dx}\Big(\f(t+x-\xi,x)\Big)\, dx\, dt\,d\xi=0\label{gondolin1}
\end{eqnarray}
by Dirichlet boundary conditions. Analogously,
\begin{equation}\label{gondolin2}
-\int_\O \bigg(
   \f_t(t,x)\int_0^\kappa  h(t-x-\xi,\xi)d\xi
   + \f_x(t,x)\int_0^\kappa  h(t-x-\xi,\xi)d\xi
   \bigg)\,dt\, dx=0\,.
\end{equation}
Moreover, again by Fubini Theorem,
\begin{eqnarray}
   & &  \int_\O \bigg(-
   \f_t(t,x)\int_\kappa^x  h(t+x-\xi,\xi)d\xi
   + \f_x(t,x)\int_\kappa^x  h(t+x-\xi,\xi)d\xi
   \bigg)\,dt\, dx \nonumber\\
   & & =
   \int_0^\pi \int_\kappa^x \intt  \Big(
   - \f_t(t,x)  h(t+x-\xi,\xi)
   + \f_x(t,x)  h(t+x-\xi,\xi)
   \Big)\,dt\,d\xi\, dx
   \nonumber\\
   & &=
   \int_0^\pi \int_\kappa^x \intt  \Big(
   - \f_t(t-x+\xi,x)+
   \f_x(t-x+\xi,x)
   \Big)h(t,\xi)\,dt\,d\xi\, dx
   \nonumber\\
   & & =\intt \int_0^\pi h(t,\xi) \int_\xi^\pi
   \frac{d}{dx}\Big(\f(t-x+\xi,x)\Big)\, dx\, dt\,d\xi\nonumber\\
   & &
   \ \ \ -\intt \int_0^\kappa h(t,\xi) \int_0^\pi
   \frac{d}{dx}\Big(\f(t-x+\xi,x)\Big)\, dx\, dt\,d\xi\nonumber
= - \int_\O h \f
   \label{gondolin3}
\end{eqnarray}
and, analogously,
\begin{equation}\label{gondolin4}
\int_\O \bigg(
   \f_t(t,x)\int_\kappa^x  h(t-x+\xi,\xi)d\xi
   + \f_x(t,x)\int_\kappa^x  h(t-x+\xi,\xi)d\xi
   \bigg)\,dt\, dx=- \int_\O h \f\,.
\end{equation}
Summing
\equ{gondolin1}, \equ{gondolin2}, \equ{gondolin3}, \equ{gondolin4}
and recalling \equ{erika}, \equ{ale} we get \equ{raffi}.
\\[1mm]
{\bf Step 4:} {\sl $ H(t,x) > 0 $ in $ \O $.}
\\[1mm]
{\sl First case: $ 0 < x \leq \kappa $.}
By  \equ{Hc1} and geometrical considerations on the
domains of the integrals, we derive that,
for $ 0 < x < \kappa $,
$ H(t,x) = \int_{\Theta} h (\tau, \xi) \ d \tau d \xi $
where $ \Theta := \Theta_{t,x} $ is the trapezoidal region in $ \O $ with
a vertex in $(\tau, \xi ) = (t,x) $ and delimited by
the straight lines
$ \tau = t - x + \xi $, $ \tau = t + x - \xi $, $ \x = \kappa $ and
$ \tau = t - x - \xi $. Since $ h > 0 $ a.e. in $\O$
we conclude that $ H (t, x) > 0 $.
\\[1mm]
{\sl Second case: $ \kappa  < x < \pi $.}
Since $ H(t+ \pi - x, \pi) = 0 $ 
we have, by \equ{Hc1},
$$
\int_0^\kappa \int_{t-x-\xi}^{t-x+\xi} h(\tau, \xi) \ d \tau d \xi =
\int_\kappa^\pi \int_{t-x + \xi}^{t-x-\xi+2 \pi}
h(\tau, \xi) \ d \tau d \xi \ .
$$
Therefore, substituting in \equ{Hc1},
we get, for $ \kappa < x < \pi $,
the expression $ H(t,x) = $ $\int_{\Theta} h (\tau, \xi) \ d \tau d \xi $
where, now, $ \Theta := \Theta_{t,x} $
is the trapezoidal region in $ \O $ with
a vertex in $(\tau, \xi ) = (t,x) $ and delimited by
the straight lines
$ \tau = t - x + \xi $,
$ \tau = t - x - \xi + 2 \pi $, $ \xi = \kappa $ and
$ \tau = t + x - \xi $. Since $ h > 0 $ a.e. in $\O$
we conclude also in this case that $ H (t, x) > 0 $.
\qed

\medskip
\noindent
\textsc{Proof of Theorem 1.}
Since $h>0$ a.e. in $\O$, by Proposition \ref{lem:6}
there exists a weak solution $ H \in E $ of
$\square H = h $ verifying \equ{G:positiva} (i.e. $ H > 0 $ in $ \O $).
Therefore existence of a weak solution $ u \in E $
satisfying $ \| u \|_E \leq C |\e | $
follows from Theorem 2-(i)
with $ \b(x) \equiv \b $ and $ \mathcal{R} \equiv 0 $.
The higher regularity for $ u $
and the estimate $\| u\|_{H^{j+1}( \O)}+$ $
\| u\|_{C^j(\oO)}\leq C |\e|$
follow from Theorem 2-(ii) and
\equ{sti:regol} in remark \ref{rem:Hp3}
since assumption \equ{hbis} is trivially verified (${\cal R} \equiv 0$).
\qed

\section{Proof of Theorem 3}\label{sec:thm2}

In order to prove Theorem 3
we perform the Lyapunov-Schmidt reduction of section \ref{sec:L-S}
and we minimize
the reduced
action functional $\Phi$
in $ \ov{B_R} := \{ \|v\|_{H^1}\leq R\}$.
To conclude the existence of a solution,
we have to prove that the minimum $ \bv \in \ov{B_R}$ is an interior
minimum point in $B_R $ for some $R > 0 $.

This case is easier that the previous one since
the required a-priori estimates can be deduced directly by the
$ 0^{th}$-order variational inequality \equ{minore}
which does not vanish for $ \e = 0 $.
\\[2mm]
{\bf Step 1: The $L^\infty$-estimate}
\\[1mm]
\indent
Since $a(x,u)$ satisfies \equ{a1} or \equ{a2},
by \equ{pari},  $ \int_\O a(x,\bv)\f=0,$
$\forall\f\in V $ (as $ \bv\in V $) and hence
\begin{equation}\label{zero}
\int_\O f(\bv)\f=\int_\O \tf(\bv)\f + a (x, \bv ) \f =
\int_\O \tf(\bv)\f \ .
\end{equation}
Since $ \|w (\bv, \e) \|_E \leq C |\e | $,
by the variational inequality
\equ{minore} and \equ{zero}, we find
$$
\int_\O \tf(\bv)\f\leq
o(1)C_1(\| \bv\|_{L^\infty})
\| \f\|_{L^1}
\leq o(1)C_1(R)\| \f\|_{L^1}\,,
$$
where $C_1(\cdot)$ is a suitable increasing function
depending on
$f.$
Then there exists a
decreasing function $0<\e_1(\cdot)\leq \e_0(\cdot)$ such that
\begin{equation}\label{stima:Li1R1}
\int_\O \tf(\bv)  \f\leq \| \f\|_{L^1}
\, \qquad\mbox{for}\quad
|\e|\leq \e_1(R)\,.
\end{equation}
We now choose, as in subsection \ref{subsec:Li}, the admissible variation
$\f = q( \bv_+ ) - q(\bv_- ) $ where $q$ is defined in \equ{q}.
By the mean value Theorem
$$
\tf(t,x,\bv)=\tf(t,x,0)+\tf_u(\mbox{intermediate point})\bv\,
$$
and, by \equ{stima:Li1R1},
since $ \tf_u \geq \b > 0 $
and $ \bv \f \geq 0 $ (recall \equ{vfpos}),
we obtain
\begin{equation}\label{stima:Li2R1}
\b\int_\O \bv  \f\leq \kappa_{1}\| \f\|_{L^1}
\,.
\end{equation}
Arguing as at the end of subsection \ref{subsec:Li} (see inequality
\equ{pipi2}, recall \equ{fq} and
$ M := \| \hat \bv \|_{L^{\infty}(\torus )} \slash 2$)
we deduce
$$
\int_\O \bv  \f\geq \kappa_{2}\| \bv \|_{L^\infty} \| \f\|_{L^1}
$$
and, by \equ{stima:Li2R1}, we deduce
$$
\| \bv\|_{L^\infty} \leq
\kappa_{3}
\, \qquad\mbox{for}\quad
|\e|\leq \e_1(R)\,.
$$
\\[1mm]
{\bf Step 2: The $H^1$-estimate}
\\[1mm]
\indent
The  $H^1$-estimate is carried out
as in  subsection \ref{subsec:LH1}
taking the admissible variation $\f:=-D_{-h}D_{h}\bv $
in the variational inequality \equ{minore}.
By Lemma \ref{lem:quoz},
denoting $\bw:=w(\bv)$,
\begin{equation}\label{stima:LH1R1}
0\geq \int_\O f (\bv+\bw) \f =
\int_\O D_h f (\bv+\bw) D_h\bv
\stackrel{h \to 0} \longrightarrow
\int_\O \Big( f_t ( \bv + \bw )
+  f_u (\bv+\bw)(\bv_t + \bw_t)\Big)\bv_t\,.
\end{equation}
Since $ \| w \|_E = O( \e ) $, we get,
by \equ{stima:LH1R1},
\begin{equation}\label{stima:LH1R1bis}
\int_\O f_u(\bv+\bw)\bv_t^2
\leq \kappa_{4} \| \bv_t\|_{L^2}\ .
\end{equation}
Since $\tf_u\geq \b > 0 $,
$ \int_\O a_u (x,\bv)\bv_t^2=0 $ by \equ{pari2}, and $\|w\|_{E}= O(\e )$
\begin{eqnarray}
\int_\O f_u(\bv+\bw)\bv_t^2 & = &
\int_\O \tf_u ( \bv + \bw ) \bv_t^2
+ a_u ( t, \bv + \bw ) \bv_t^2  \nonumber \\
& \geq & \int_\O \b \bv_t^2 +
\int_\O \Big( a_u ( t, \bv + \bw ) - a_u ( t, \bv ) \Big) \bv_t^2 \nonumber \\
& \geq & \int_\O \b \bv_t^2 - 
o(1) \int_\O \bv_t^2 \geq \frac{\b}{2} \int_\O \bv_t^2 \label{dalba}
\end{eqnarray}
for $ | \e | \leq \e_2(R) \leq \e_1 (R) $. 
From \equ{stima:LH1R1bis} and \equ{dalba} we deduce
\begin{equation*}
    \|\bv\|_{H^1}<\kappa_{5} \,  \qquad \forall\, |\e|\leq \e_2 (R)\, .
\end{equation*}

\noindent
\textsc{{ Proof \ of \ Theorem \ 3 completed.}}
For $ R_* := \kappa_{5} $ and $ \e_* := \e_2 ( R_* ) $,
$\bv $ is an interior point in $ B_{R_*} $
and $ \bu := \bv + w(\bv, \e) $
is a weak solution  of \equ{equ}-\equ{bc}-\equ{pc}.
Regularity of the solution $ \bu $
is proved as in subsection \ref{sec:hig}
and Theorem 3 follows.
\qed

\section{Appendix}

\noindent
\textsc{Proof \ of \ Lemma \ref{lem:0}}.
By the periodicity of $ a ( t, x ) $ with respect to $ t $
$$
\int_{\O_\a} a(t,x) dt dx = \int_{\widetilde \O_\a}  a(t,x) dt dx
$$
where $\widetilde \O_\a :=
\{ \a \pi < x < \pi (1- \a ), - x < t < -x + 2 \pi  \}$.
Under the change of variables $ s_+ := t + x  $, $ s_- := t - x  $ the domain
$\widetilde \O_\a$ trasforms into the domain
$$
\Big\{  0 < s_+ < 2 \pi, \   s_+ - 2 \pi (1 -\a ) < s_- < s_+ - 2 \pi \a  \Big\}
$$
and we get \equ{fubini1}.

For $ p, q \in L^1(\mathbb{T}) $, by \equ{fubini1} we have
\begin{eqnarray*}
\intOa p(t+x) q(t-x)\,dt\,dx & = & \frac{1}{2} \intt ds_+\ p( s_+ )
\int_{-2\pi+s_++2\a\pi}^{s_+-2\a\pi} q( s_- ) ds_- \\
& = & \frac{1}{2}  \intt ds_+\ p( s_+ ) \ \Big( \intt q( s ) ds \ -
\int_{s_+-2\a\pi}^{s_+ + 2\a\pi} q( s_- ) ds_- \Big) \\
& = & \frac{1}{2} \intt p( s ) ds \ \intt q( s ) ds \\
& & -
\frac{1}{2} \intt ds_+\ p( s_+\ )
\int_{-2\a\pi}^{2\a\pi} q( s_+ + z ) dz
\end{eqnarray*}
and we obtain \equ{fubini2} by Tonelli's Theorem (calling $s_+ = y$).

\noindent
Formula \equ{fubini6} follows by \equ{fubini2} setting $q \equiv 1$ .

We now prove \equ{fubini7}.
Since the change of variables
$(t,x) \mapsto (t,\pi-x) $ leaves
the domain $\O_\a $ unchanged
$$
    \intOa a(t,x)dtdx=
    \intOa a(t,\pi-x)dtdx\ ,
$$
and, using also the periodicity of $ p $,
\begin{eqnarray*}
  \int_{\O_\a} f(p(t+x)) g(p(t-x))\,dt\,dx  &=&
  \int_{\O_\a} f(p(t+\pi-x)) g(p(t-\pi+x))\,dt\,dx  \\
&=& \int_{\a\pi}^{\pi-\a\pi} \intt
   f(p(t+\pi-x)) g(p(t + \pi+x))\,dt\,dx  \\
  &=&
  \int_{\a\pi}^{\pi-\a\pi} \intt f(p(t-x)) g(p(t+x))\,dt\,dx
\end{eqnarray*}
proving \equ{fubini7}. \qed

\medskip\medskip

\noindent
\textsc{Proof \ of \ Lemma \ref{iso}}.
\equ{isometria} follows from the equality
$v^2=$ $v_+^2+v_-^2-2v_+v_-$,
\equ{fubini6} (with $p=v_+^2 ,v_-^2 $ and $\a=0$) and \equ{fubini2}
(taking $p,q=\hv$, $\a = 0 $ and  recalling that $\hv$ has zero average).

\equ{isometriaH} follows form \equ{isometria}
since $ v_t(t,x) = \hv' (t+x) - \hv' (t-x)$ (and similarly for $ v_x $).

Next, the first inequality of \equ{isometriaL} follows from
$ v (t,x) = \hv (t+x) - \hv (t-x)$ recalling
that, since $\hv$ has zero average,
there exist two positive measure sets in which $\hv\geq 0$ and $\hv \leq 0.$
The second inequality of \equ{isometriaL} is trivial.

We finally prove  \equ{embed}.
Since $\hv$ is continuous ($\hv\in H^1(\mathbb{T})$)
there exists $\xi_M$  such that
$\| \hv\|_{L^\infty(\mathbb{T})}$ $=|\hv(\xi_M)|$.
Being $\intt \hv=0$, there exists $|\xi_0-\xi_M|\leq \pi$
such that $\hv(\xi_0)=0$.
Hence
$$
\| \hv\|_{L^\infty(\mathbb{T})}=|\hv(\xi_M)|
= \left| \int_{\xi_0}^{\xi_M} \hv'(s)ds  \right|
\leq \sqrt{\pi} \| \hv'\|_{L^2(\mathbb{T})}
= \frac{1}{\sqrt 2} \| v_t\|_{L^2(\O)}
$$
by the Cauchy-Schwartz inequality and \equ{isometriaH}.
Finally by \equ{isometriaL},
$$
\|v\|_{L^\infty ( \O )} \leq 2 \|\hv\|_{L^\infty ( \mathbb{T})}
\leq \sqrt{2} \| v_t\|_{L^2(\O)}
= \frac{1}{\sqrt 2} \Big(\| v_t\|_{L^2(\O)} +  \| v_x\|_{L^2(\O)} \Big)
\leq \| v\|_{H^1(\O)}\,.
$$
where $\| v\|_{H^1(\O)}^2 := $
$\| v \|_{L^2(\O)}^2 + $ $\|v_x\|_{L^2(\O)}^2 + $ $ \|v_t\|_{L^2(\O)}^2 $.
\qed

\medskip\medskip

\noindent
\textsc{Proof \ of \ Lemma \ref{lem:2}}.
By the change of variables $(t,x) \mapsto (t,\pi-x) $
and periodicity,
\begin{eqnarray*}
  \intOa \f_1\cdot\ldots\cdot \f_{2k+1} &=&
  \int_{\a\pi}^{\pi-\a\pi} \intt \prod_{j=1}^{2k+1}
  \Big(  \hf_j(t+x)-\hf_j(t-x) \Big)dt dx \\
   &=&  \int_{\a\pi}^{\pi-\a\pi} \intt \prod_{j=1}^{2k+1}
  \Big(  \hf_j(t+\pi -x)-\hf_j(t-\pi+x) \Big)dt dx \\
   &=& \int_{\a\pi}^{\pi-\a\pi} \intt \prod_{j=1}^{2k+1}
  \Big(  \hf_j(t-x)-\hf_j(t+x) \Big)dt dx  \\
   &=& (-1)^{2k+1} \intOa \f_1\cdot\ldots\cdot \f_{2k+1}\ ,
\end{eqnarray*}
which implies \equ{fubini4}.

\noindent
With similar arguments we can prove
\equ{pari} and \equ{pari2}. \qed

\medskip\medskip

\noindent
\textsc{Proof \ of \ Lemma \ref{lem:ineq}}.
The inequality \equ{dis:ovvia} follows by the convexity of $t \to t^{2k}$.

We next prove \equ{zialalletta}.
If $b=0$ it is trivially true. If $b\neq 0$ let us divide for $b^{2k}$ and
set $ x:= a/b \in\mathbb{R}$. \equ{zialalletta} is equivalent to prove
\begin{equation}\label{disfu}
f(x):=(x-1)^{2k}-x^{2k}-1+2k x^{2k-1}+2kx \geq 0.
\end{equation}
It results
$f(0)=0$, $f'(x)=2k[(x-1)^{2k-1}-x^{2k-1}+(2k-1)x^{2k-2}+1]$ and
so $f'(0)=0$. Therefore to prove \equ{disfu}
it is sufficient to show that $f$ is convex. We have
$f''(x)=2k(2k-1)g(x)$ where
$ g_{k}(x) := (x-1)^{2k-2}-x^{2k-2}+(2k-2)x^{2k-3}$, $k\geq2$.\\
We now show by induction on $k\geq2$ that $g_k (x) > 0$. It is true for
$k=2$ since $ g_{2}(x)= $ $(x-1)^{2}-x^{2}+2x=$ $ 1 > 0.$
Supposing now $ g_k (x) > 0 $, let us prove that $g_{k+1}(x) > 0 $.

We claim that
$$g_{k+1}(x)=(x-1)^{2k}-x^{2k}+2kx^{2k-1} $$
is a strictly convex function. Indeed
$$ g_{k+1}'(x)=2k[(x-1)^{2k-1}-x^{2k-1}+(2k-1)x^{2k-2}] \quad
{\rm and} \quad
g_{k+1}''(x)=2k(2k-1)g_{k}(x).
$$
By the inductive hypothesis $g_{k}(x)>0$ and therefore
$g_{k+1}''(x)>0$. Moreover, being
$g_{k+1}(x)\approx \mathrm{cost} x^{2k-2}$,
$ \lim_{x\rightarrow \pm \infty}g_{k+1}(x)=$ $+\infty $
and $g_{k+1}(x)$ possesses a unique point of global minimum $\bar{x}$ that is
also the unique critical point. Now it is sufficient to show that
$g_{k+1}(\bar{x})>0$.
$$g_{k}'(\bar{x})=2k[(\bar{x}-1)^{2k-1}-\bar{x}^{2k-1}+(2k-1)\bar{x}^{2k-2}] =0$$
implies that $(\bar{x}-1)^{2k-1}= $ $\bar{x}^{2k-1}-(2k-1)\bar{x}^{2k-2}.$
Substituting this equality in the expression for $g_{k+1}(\bar{x})$,
we have
\begin{eqnarray*}
  g_{k+1}(\bar{x}) &=& (\bar{x}-1)[\bar{x}^{2k-1}-(2k-1)\bar{x}^{2k-2}]-
  \bar{x}^{2k}+2k\bar{x}^{2k-1}  \\
   &=& \bar{x}^{2k}-\bar{x}^{2k-1}-(2k-1)\bar{x}^{2k-1}+(2k-1)\bar{x}^{2k-2}-
  \bar{x}^{2k}+2k\bar{x}^{2k-1}  \\
  &=& (2k-1)\bar{x}^{2k-2}>0
\end{eqnarray*}
(we use that $\bar{x}\neq 0$, in fact $g_{k+1}'(0)=-2k\neq 0)$).
\\[1mm]
\indent
Proof of \equ{zialalletta3}.
The case $k=1$ is trivial.
For $k\geq 2$, we divide by $b^{2k-1}$ and define $x:=a/b\in\mathbb{R}.$
We have to prove that
$$
f(x):=(x+1)^{2k-1}-x^{2k-1}\geq 4^{1-k}\ ,\qquad\ \ \ \forall\ x\in\mathbb{R}\ .
$$
Since
$$
f'(x)=(2k-1)[(x+1)^{2k-2}-x^{2k-2}]=0\quad \Longleftrightarrow \quad
(x+1)^{2k-2}=x^{2k-2}\quad \Longleftrightarrow \quad
x=-\frac{1}{2}
$$
and $f(x)\to \infty$ as $|x|\to\infty$,
we conclude that $x=-1/2$ is the unique minimum point of $f(x)$ and
therefore $f(x)\geq f(-1/2)=4^{1-k}.$ \qed

\medskip\medskip

\noindent
\textsc{Proof of Lemma \ref{lem:3}}.
From the inequality \equ{dis:ovvia} we obtain
$$
\int_\O v^{2k} =
\int_\O (v_+ - v_-)^{2k}
\leq 2^{2k-1}\int_\O v_+^{2k}+v_-^{2k}
=
2^{2k-1}\int_\O \hv^{2k}(t+x)+\hv^{2k}(t-x)dtdx\ ,
$$
which, using
  \equ{fubini6},
proves \equ{fubini5}.\qed

We first prove the Lemma in the case $k\geq 2.$
Using the inequality \equ{zialalletta}
\equ{fubini6} and \equ{fubini7}, we obtain
\begin{eqnarray}
\intOa v^{2k} &=&
\intOa (v_+ -v_-)^{2k}\geq
\intOa v_+^{2k}  +v_-^{2k} -2k \intOa v_+^{2k-1} v_-
+ v_+ v_-^{2k-1}\nonumber\\
&=&
 2\pi (1-2\a)\intt \hv^{2k} -4k \intOa v_+^{2k-1} v_-\ .
\label{zialalla}
\end{eqnarray}
By \equ{fubini2} and since $\hv $ has zero average
\begin{eqnarray}\label{zialallissima}
\left|\intOa v_+^{2k-1} v_-\right| &=&
\frac12\left|
\int_{-2\pi\a}^{2\pi\a}
\intt \hv^{2k-1}(y)\hv(y+z)dydz\right|\nonumber\\
&\leq&
\frac12
\int_{-2\pi\a}^{2\pi\a}
\intt |\hv(y)|^{2k-1}|\hv(y+z)|dydz\ .
\end{eqnarray}
By H\"older inequality with $p:=2k/(2k-1)$ and $q:=2k$ ($1/p+1/q=1$),
\begin{eqnarray*}
\intt |\hv(y)|^{2k-1}|\hv(y+z)|dy &\leq&
\left(\intt |\hv(y)|^{2k}dy\right)^{\frac{2k-1}{2k}}
\left(\intt |\hv(y+z)|^{2k}dy\right)^{\frac{1}{2k}}\\
&=&
\intt |\hv(y)|^{2k}dy\ ,
\end{eqnarray*}
where, in the equality, we have used the periodicity of $\hv$
to conclude that $\intt |\hv(y+z)|^{2k}dy$
$=\intt |\hv(y)|^{2k}dy$.
Inserting the last inequality in \equ{zialallissima},
we obtain
\begin{equation}\label{lalla}
\left|\intOa v_+^{2k-1} v_-\right|
\leq
2\pi\a \intt |\hv(y)|^{2k}dy
\end{equation}
Inserting \equ{lalla} in \equ{zialalla} the Proposition follows
in the case $k\geq 2.$
 The case $k=1$ is similar:
$$
\intOa v^2 = \intOa v_+^2+v_-^2 -2v_+ v_- =
 2\pi (1-2\a)\intt \hv^2 -2 \intOa v_+ v_-
$$
and we conclude by \equ{lalla}.\qed

\medskip\medskip

\noindent
\textsc{Proof of Lemma \ref{lem:quoz}}.
Formula \equ{Leibniz} follows from
\begin{eqnarray*}
D_h(fg)(t,x) &= & \frac{f(t+h,x)g(t+h,x)-f(t,x)g(t,x)}{h} \\
& = & \frac{f(t+h,x)-f(t,x)}{h}g(t,x)+
f(t+h,x)\frac{g(t+h,x)-g(t,x)}{h}\ .
\end{eqnarray*}
We prove \equ{Leibniz2} by induction. It is obvious for $m=1$.
We suppose \equ{Leibniz2} holds for $m$ and prove it for $m+1:$
by \equ{Leibniz} we have
\begin{eqnarray*}
D_h (f^m \cdot f) &=&
         (D_h f^m) f+ T_h f^m D_h f \\
&=&
(D_h f) \left[
\sum_{j=0}^{m-1} f^{m-j} T_h f^j + T_h f^m
 \right]\\
&=&
(D_h f)
\sum_{j=0}^{m} f^{(m+1)-j-1} T_h f^j \ .
\end{eqnarray*}
Formula \equ{Leibnizint} follows by \equ{Leibniz}
performing the change of variables $s=t+h$,
\begin{eqnarray*}
  \int_\O D_h(fg)(t,x) &=&
  \int_\O (D_h f) g +
\int_\O f(t+h,x) \frac{g(t+h,x)-g(t,x)}{h} dtdx \\
   &=& \int_\O (D_h f) g +
\int_\O f(s,x) \frac{g(s,x)-g(s-h,x)}{h} dsdx\,.
\end{eqnarray*}
We now prove formula \equ{perparti} for integration by parts.
Due to the periodicity of $f$ and $g$ with respect to  $ t $
\begin{eqnarray*}
\int_{\O}
f (D_{-h}g)  &=&
-\frac{1}{h}\int_{\O}
 f(t,x) \Big[ g(t-h,x)-g(t,x) \Big] dtdx\\
&=& -\frac{1}{h}
\int_{\O} \Big[ f(t+h,x)-f(t,x) \Big] g(t,x)dtdx
=
-
\int_{\O}
(D_{h}f) g\ .
\end{eqnarray*}

The proof of \equ{weakd} is standard.
Let $ \f \in C^1 ( \O ). $ By \equ{perparti}
\begin{equation}\label{perpa}
\int_\O (D_h f)\f=-\int_\O f (D_{-h}\f)\ .
\end{equation}
Now the sequence $f (D_{-h}\f)$ converges to $ f\f_t$ a.e. and,
since
$$
|f(t,x)\, (D_{-h}\f)(t,x)|\leq \| \f \|_{C^1(\O)} |f(t,x)|
\in L^1(\O) \ ,
$$
we can use the Lebesgue Theorem to obtain
\begin{equation}\label{con}
\int_\O f (D_{-h}\f) \stackrel{h \to 0 }\longrightarrow
\int_\O f \f_t.
\end{equation}
Since, by hypothesis, $D_h f$ is bounded in $L^2 (\O ) $,
$ D_h f \stackrel{L^2}\rightharpoonup g $, up to a subsequence.
Passing to the limit in \equ{perpa} for $ h \to 0 $
we find
$ \int_\O g \f= $ $ - \int_\O f \f_t.$ Therefore $f$
has a weak derivative $f_t = g $ and
by the weakly lower semicontinuity of the norm
$$
\|f_t \|_{L^2} \leq \liminf \| D_h f \|_{L^2} \leq C \ .
$$

In order to prove \equ{sti:diffquo} assume temporarily $f$ is smooth.
From the fundamental Theorem of calculus
$$
(D_h f)(t,x) =\frac{f(t+h,x)-f(t,x)}{h}=\int_0^1 f_t(t+hs,x) ds\ .
$$
By Cauchy-Schwartz inequality, Fubini Theorem and periodicity
we obtain
\begin{eqnarray*}
\int_\O |D_h f(t,x)|^2 dtdx&=& \int_0^\pi \intt \Big|
\int_0^1 f_t(t+hs,x) ds
 \Big|^2 dtdx \\
&\leq&
\int_0^\pi \intt
\int_0^1 | f_t(t+hs,x)
 |^2 dsdtdx \\
&=&
\int_0^\pi
\int_0^1
\intt
|
f_t(t+hs,x)
 |^2 dt ds dx\\
&=&
\int_0^\pi
\int_0^1
\intt
|
f_t(t,x)
 |^2 dt ds dx\\
&=&
\int_0^1
\int_0^\pi
\intt
|
f_t(t,x)
 |^2 dt dx ds
=
\| f_t \|_{L^2(\O)}^2 \ .
\end{eqnarray*}
Inequality \equ{sti:diffquo} is valid,
for any $f$ having a weak derivative $ f_t \in L^2 (\O ) $,
by approximation.

In order to prove \equ{conL2quoz}
we first show the weak $L^2$-convergence.
Let $\f\in C^1(\O).$
By \equ{perparti}, applying as before the
Lebesgue Theorem, and since $f $ has a weak derivative $ f_t  $
\begin{equation}\label{condeb}
\int_\O (D_h f)\f = -\int_\O f (D_{-h}\f)
\stackrel{h \to 0 }\longrightarrow -\int_\O f \f_t =
\int_\O f_t \f\ , \qquad \forall\,\f\in C^1(\O)\ .
\end{equation}
Since, by \equ{sti:diffquo}, $D_h f$ is bounded in $L^2 $, and
\equ{condeb} holds in the dense subset $C^1(\O)$ $\subset L^2(\O)$,
we conclude the weak $L^2$-convergence
$ D_h f \stackrel{L^2}\rightharpoonup $ $ f_t$.

By the weakly lower semicontinuity of the norm
and \equ{sti:diffquo}
\begin{equation}\label{convinnorm}
\| f_t\|_{L^2(\O)}\leq \liminf\|D_h f \|_{L^2(\O)}
\leq \| f_t\|_{L^2(\O)}\quad\Longrightarrow\quad
\lim \|D_h f \|_{L^2(\O)}
= \| f_t\|_{L^2(\O)}\ .
\end{equation}
Since $ L^2(\O)$ is a Hilbert space, weak convergence
$ D_h f \stackrel{L^2}\rightharpoonup $ $ f_t$ and
\equ{convinnorm} imply the strong convergence
$ D_h f \stackrel{L^2}\rightarrow $ $ f_t$.\qed

\section*{References}

\footnotesize

\noindent
[BP01] Bambusi, D.; Paleari, S. :
{\sl Families of periodic solutions of resonant PDEs},
J. Nonlinear Sci.  11  (2001),  no. 1, 69--87.
\smallskip

\noindent
[BDL99] Bartsch, T.; Ding, Y. H.; Lee, C. :
{\sl Periodic solutions of a wave equation with concave
and convex nonlinearities},  J. Differential Equations  153  (1999),
no. 1, 121--141.
\smallskip

\noindent
[BB03] Berti, M.; Bolle, P.:
{\sl Periodic solutions of nonlinear wave equations with general
nonlinearities}, Comm. Math. Phys.  243  (2003),  no. 2, 315--328.
\smallskip

\noindent
[BB04a] Berti, M.; Bolle, P.:
{\sl Multiplicity of periodic solutions of nonlinear
wave equations}, Nonlinear Anal.  56,  (2004), 1011--1046.
\smallskip

\noindent
[BB04b] Berti, M.; Bolle, P.:
{\sl Cantor families of periodic solution for completely
resonant wave equations}, preprint SISSA 2004.
\smallskip

\noindent
[B99] J. Bourgain, {\it Periodic solutions of nonlinear wave
equations},  Harmonic analysis and partial differential equations,
69--97, Chicago Lectures in Math., Univ. Chicago Press, 1999.
\smallskip

\noindent
[BN78] Br\'ezis, H.; Nirenberg, L. :
{\sl Forced vibrations for a nonlinear wave equation},
Comm. Pure Appl. Math.  31  (1978), no. 1, 1--30.
\smallskip

\noindent
[BCN80]
Br\'ezis, H.; Coron, J.-M.; Nirenberg, L.
{\sl Free vibrations for a nonlinear wave equation and a
Theorem of P. Rabinowitz},
Comm. Pure Appl. Math.  33  (1980), no. 5, 667--684.
\smallskip

\noindent
[C83] Coron, J.-M. :
{\sl Periodic solutions of a nonlinear wave equation without
assumption of monotonicity}, Math. Ann. 262 (1983), no. 2, 273--285.
\smallskip

\noindent
[DST68] De Simon, L.; Torelli, H. :
{\sl  Soluzioni periodiche di equazioni a derivate parziali
di tipo iperbolico non lineari},
Rend. Sem. Mat. Univ. Padova  40  1968 380--401.
\smallskip

\noindent
[GMP04] G. Gentile, V. Mastropietro, M. Procesi, {\sl
Periodic solutions for completely resonant nonlinear wave equations},
to appear on Comm. Math. Phys.
\smallskip

\noindent
[H70]
Hall, W. S. :
{\sl On the existence of periodic solutions for the equations
$D_{tt}u+(-1)\sp{p}\,D_{x}\sp{2p}u=\varepsilon f( \cdot, \cdot ,u)$},
J. Differential Equations,  7,  1970, 509--526.
\smallskip

\noindent
[H82]
Hofer, H.:
{\sl On the range of a wave operator with nonmonotone nonlinearity},
Math. Nachr. 106 (1982), 327--340.
\smallskip

\noindent
[L69] Lovicarov\'a, H.:
{\sl Periodic solutions of a weakly nonlinear wave equation
in one dimension}, Czechoslovak Math. J. 19 (94) 1969 324--342.
\smallskip

\noindent
[PY89]
Plotnikov, P. I.; Yungerman, L. N.:
{\sl Periodic solutions of a weakly nonlinear wave equation
with an irrational relation of period to interval length}, (Russian)
Differentsial nye Uravneniya 24 (1988), no. 9, 1599--1607, 1654;
translation in Differential Equations 24 (1988), no. 9, 1059--1065 (1989).
\smallskip

\noindent
[R67] Rabinowitz, P.:
{\sl Periodic solutions of nonlinear hyperbolic partial differential
equations}, Comm. Pure Appl. Math.  20, 145--205, 1967.
\smallskip

\noindent
[R71] Rabinowitz, P.:
{\sl Time periodic solutions of nonlinear wave equations},
Manuscripta Math.  5  (1971), 165--194.
\smallskip

\noindent
[T69]
Torelli, G.:
{\sl Soluzioni periodiche dell'equazione non lineare
$u_{tt}-u_{xx} + \varepsilon F(x,\,t,\,u)=0$},
Rend. Ist. Mat. Univ. Trieste  1  1969 123--137.
\smallskip

\noindent
[W81] Willem, M.
{\sl Density of the range of potential operators},
Proc. Amer. Math. Soc. 83 (1981), no. 2, 341--344.
\smallskip

\end{document}